\documentclass[12pt,a4paper]{article}
\usepackage[utf8]{inputenc}
\usepackage{xifthen}
\usepackage{amssymb}
\usepackage{amsthm}
\usepackage{mathtools}
\usepackage{physics} 
\usepackage{enumerate}
\usepackage{hyperref}
\usepackage{cleveref}
\usepackage{nameref}
\usepackage{comment}
\usepackage{caption} 
\usepackage{subcaption}
\usepackage{tikz}

\usepackage[TS1,T1]{fontenc}
\usepackage{libertine}
\usepackage{textcomp} 
\usepackage[varqu,varl]{inconsolata}
\usepackage[libertine,cmintegrals,varbb,slantedGreek]{newtxmath}
\usepackage[scr=rsfso]{mathalfa}
\usepackage{bm}
\useosf 
\usepackage{enumerate}
\usepackage[supstfm=libertinesups,supscaled=1.2,raised=-.13em]{superiors}
\usepackage{fnpct}
\usepackage[final]{microtype}

\usepackage[american,main=american]{babel}
\AtBeginDocument{\selectlanguage{american}}

\usepackage[autostyle,english=american]{csquotes}

\usepackage{setspace}

\newcommand\nref[1]{\nameref{#1} \ref{#1}}
\DeclareMathAlphabet\mathbfcal{OMS}{cmsy}{b}{n}

\DeclarePairedDelimiterX\set[1]\lbrace\rbrace{\let\given\undefined\newcommand\given{\;\delimsize\vert\;}#1}
\DeclarePairedDelimiter{\prn}{\lparen}{\rparen}

\newcommand{\tidx}[1]{{(#1)}}

\newcommand{\N}{\mathbb{N}}
\newcommand{\R}{\mathbb{R}}
\newcommand{\collc}{\mathbfcal{C}}
\newcommand{\collp}{\mathcal{P}}

\newcommand{\mH}[1]{\mathcal{H}^{#1}}
\DeclareMathOperator{\im}{im}
\newcommand{\diff}{\mathop{}\!\mathrm{d}}

\newcommand{\ball}[1]{{B^{#1}}}
\newcommand{\stdb}{\ball n}

\newcommand{\sphere}[1]{\mathbb{S}^{#1}}
\newcommand{\stdsph}{\sphere{n-1}}

\newcommand{\mres}[1]{\llcorner #1}

\DeclareMathOperator{\relint}{relint}

\DeclareMathOperator{\absint}{int}
\DeclareMathOperator{\cl}{cl}
\DeclareMathOperator{\V}{V}
\DeclareMathOperator{\Su}{S}

\DeclareMathOperator{\conv}{conv}

\DeclareMathOperator{\diam}{diam}
\DeclareMathOperator{\vspan}{span}
\DeclareMathOperator{\pspan}{\overline{span}}
\DeclareMathOperator{\supp}{supp}
\DeclareMathOperator{\ext}{ext}

\DeclareMathOperator{\TS}{TS}
\DeclareMathOperator{\strip}{prune}
\DeclareMathOperator{\prune}{prune}
\newcommand{\elminus}[1]{\setminus\set{#1}}
\newcommand{\ccone}{\mathfrak{C}}

\newcommand{\K}{\mathcal{K}}
\newcommand{\Pt}{\mathcal{P}}
\newcommand{\mbS}{\mathbb{S}}

\newcommand{\Qseq}{\mathbf{Q}}

\newcommand{\Vseq}{\mathbf{V}}

\newtheorem{baseenv}{Abstract Env}[section]
\newtheorem{theorem}[baseenv]{Theorem}
\newtheorem{definition}[baseenv]{Definition}
\newtheorem{remark}[baseenv]{Remark}
\newtheorem{lemma}[baseenv]{Lemma}
\newtheorem{corollary}[baseenv]{Corollary}
\newtheorem{example}[baseenv]{Example}
\newtheorem{proposition}[baseenv]{Proposition}
\newtheorem*{theorem*}{Theorem}

\newcommand{\csupp}{\supp\Su(\collc)}
\newcommand{\cextdirs}{\ext\collc}
\newcommand{\summandcolls}{{\collp\in\prod_{i=1}^{n-1}\supp\mu_i}}

\renewcommand{\complement}{\mathsf{c}}

\title{The support of mixed area measures involving a new class of convex bodies}
\author{Daniel Hug and Paul A. Reichert}
\date{\today}

\begin{document}

\maketitle

\begin{abstract}
    Mixed volumes in $n$-dimensional Euclidean space are  functionals of $n$-tuples of convex bodies $K,L,C_1,\ldots,C_{n-2}$. The Alexandrov--Fenchel inequalities are fundamental inequalities between mixed volumes of convex bodies. As very special cases they cover or imply many important inequalities between basic geometric functionals. A complete characterization of the equality cases in the Alexandrov--Fenchel inequality remains a challenging open problem. Major recent progress was made by Yair Shenfeld and Ramon van Handel   \cite{SvH22,SvH23+}, in particular they resolved the problem in the cases where $C_1,\ldots,C_{n-2}$ are polytopes, zonoids or smooth bodies (under some dimensional restriction).
    In \cite{HugReichert23+} we introduced the class of polyoids, which are defined as limits of finite Minkowski sums of polytopes having a bounded number vertices. Polyoids encompass polytopes, zonoids and  triangle bodies, and they can be characterized by means of generating measures. Based on this characterization and Shenfeld and van Handel's contribution, we extended their result to  polyoids (or  smooth bodies). Our previous result was stated in terms of the support of the mixed area measure associated with the unit ball $B^n$ and $C_1,\ldots,C_{n-2}$. This characterization result is completed in the present work which more generally provides a geometric description of the support of the mixed area measure of an arbitrary $(n-1)$-tuple of polyoids (or smooth bodies). The result confirms a long-standing conjecture by Rolf Schneider in the case of polyoids and hence, in particular, of zonoids.
\end{abstract}

\paragraph{\small MSC-classes 2020.}{\small 52A39, 52A20, 52A21, 52A40}
\paragraph{\small Keywords.}{\small Polytope, zonoid, polyoid, Alexandrov--Fenchel inequality, generating measure, mixed area measure}

\section{Introduction}\label{sec:1}

Mixed volumes of convex bodies (nonempty compact convex sets) in Euclidean space $\R^n$, $n\ge 2$, are symmetric functionals of $n$-tuples of convex bodies, which naturally arise
as coefficients of polynomial expansions of nonnegative Minkowski combinations of convex bodies. We write $\V$ for the volume functional (Lebesgue measure) and $\alpha_1K_1+\cdots+\alpha_mK_m$ for the Minkowski combination of the convex bodies $K_1,\ldots,K_m\subset\R^n$ with nonnegative coefficients $\alpha_1,\ldots,\alpha_m\in\R$. Then
\begin{equation}\label{eq:1.1}
\V(\alpha_1K_1+\cdots+\alpha_mK_m)=\sum_{i_1,\ldots,i_n=1}^m \V(K_{i_1},\ldots,K_{i_n})\alpha_{i_1}\cdots\alpha_{i_n},
\end{equation}
where $\V(K_{i_1},\ldots,K_{i_n})$ is called the mixed volume of $K_{i_1},\ldots,K_{i_n}$. A local counterpart of the mixed volumes are the mixed area measures. For convex bodies $K_1,\ldots,K_{n-1}\subset\R^n$, the mixed area measure $\Su(K_1,\ldots,K_{n-1},\cdot)$ is the uniquely determined Borel measure on the Euclidean unit sphere $\mbS^{n-1}$ such that
\begin{equation}\label{eqex}
\V(K_1,\ldots,K_{n-1},K_n)=\frac{1}{n}\int_{\mathbb{S}^{n-1}} h_{K_n}(u)\, \Su(K_1,\ldots,K_{n-1},\diff u)
\end{equation}
holds for all convex bodies $K_n\subset \R^n$, where $h_{K_n}$ is the support function of $K_n$ (see \cite[Sect.~5.1]{Schneider} or \cite[Thm.~4.1]{Hug}).

A deep inequality for mixed volumes of convex bodies, with many consequences and applications to diverse fields, has been found and established by Alexandrov \cite{AF1937} (see Schneider \cite[Sect. 7.3]{Schneider}, also for some historical comments). We write $\K^n$ for the set of convex bodies in $\R^n$.

\begin{theorem*}[Alexandrov--Fenchel Inequality]
    \makeatletter\def\currentlabelname{Alexandrov--Fenchel Inequality}\makeatother\label{thm:af}
    Let  $K, L\in\K^n$ be convex bodies, and let $\collc = (C_1, \ldots, C_{n-2})$ be an $(n-2)$-tuple of convex bodies in $\R^n$. Then
    \begin{align}
      \V(K, L, \collc)^2 \ge \V(K, K, \collc) \V(L, L, \collc),\tag{AFI}
    \end{align}
where $\V(K, L, \collc):=\V(K,L,C_1, \ldots, C_{n-2})$.
\end{theorem*}

While the inequality was already established by Alexandrov and various proofs of the inequality are known, some of which were found recently (see \cite{Cor19,SvH19,Wa18} and the references given there), a complete characterization of the equality cases remains a major open problem in Brunn--Minkowski theory (see \cite[Sect.~7.6]{Schneider}). For recent progress, we mention the work by Shenfeld and van Handel  \cite{SvH22,SvH23+} and the literature cited there. Based on their findings for the case where $\collc = (C_1, \ldots, C_{n-2})$ is a tuple of polytopes, zonoids or smooth bodies (satisfying a weak dimensionality assumption, called supercriticality), the following more general result has been shown in \cite{HugReichert23+}.
It confirms a conjecture by Rolf Schneider \cite[Conjecture 7.6.16]{Schneider} for a new class of convex bodies, which we called polyoids, that contains all polytopes, zonoids and triangle bodies. A \emph{polyoid} is a convex body $K$ for which there is some integer $k\in\N$ and a sequence of Minkowski sums of polytopes each having at most $k$ vertices that converges to $K$; see \cite[Sect.~2]{HugReichert23+} (and Section \ref{sec:3} below) for further details and a representation theorem characterizing polyoids. A convex body is smooth if each of its boundary points is contained in a unique supporting hyperplane.

\begin{theorem*}[Equality cases in (AFI) for polyoids and smooth  bodies \cite{HugReichert23+}]\label{corfin}
Let $K,L\in \K^n$,  and let $\collc = (C_1, \ldots, C_{n-2})$ be a supercritical $(n-2)$-tuple of polyoids or smooth convex bodies in $\R^n$. Assume that $\V(K,L,\collc)>0$. Then {\rm (AFI)} holds with equality if and only if there are $a>0$ and $x\in\R^n$ such that
$$h_K=h_{aL+x}\quad\text{ on }\supp \Su(\stdb,\collc,\cdot),$$
where $\supp \Su(\stdb,\collc,\cdot)$ denotes the support of the mixed area measure $\Su(\stdb,\collc,\cdot)$ of the unit ball $B^n$ and the $(n-2)$-tuple $\collc$.
\end{theorem*}

For a geometric understanding of the equality cases in the Alexandrov--Fenchel inequality (AFI) it thus remains to describe the support of the measure $\Su(\stdb,\collc,\cdot)$ in geometric terms. According to another (more general) conjecture by Rolf Schneider \cite[Conjecture 7.6.14]{Schneider}, the support of the mixed area measure $\Su(K_1,\ldots,K_{n-1},\cdot)$, for given convex bodies $K_1,\ldots,K_{n-1}\subset\R^n$, is the closure of the set of \emph{$(K_1,\ldots,K_{n-1})$-extreme normal vectors}, for which we write $\cl \ext(K_1,\ldots,K_{n-1})$; an explicit definition and further information are given in Section \ref{sec:2}. If all convex bodies are polytopes or all are smooth and strictly convex, then the conjecture is known to be true. The conjecture was also
recently confirmed by Shenfeld and van Handel in the case of $(n-1)$-tuples of the form $(B^n, C_1,\ldots,C_{n-2})$, where $C_i$ is a zonoid or a smooth convex body in $\R^n$. However, even in the case where the unit ball $B^n$ is replaced by a general zonoid, the conjecture was open up to now.

Our main result confirms Schneider's conjecture \cite[Conjecture 7.6.14]{Schneider} not only for general $(n-1)$-tuples of zonoids (or smooth bodies), but for the larger class of polyoids (or smooth bodies).

\begin{theorem}
  \makeatletter\def\@currentlabelname{Support Characterization Theorem}\makeatother\label{thm:suppChar}
  Let $\collc = (C_1, \ldots, C_{n-1})$ be an $(n-1)$-tuple of polyoids (or smooth convex bodies provided at least one of the bodies $C_i$ is smooth and strictly convex) in $\R^n$. Then
  \begin{equation}\label{eq:key}
 \supp\Su(\collc,\cdot) = \cl\cextdirs
  .\end{equation}
\end{theorem}

In combination with the preceding theorem on the characterization of the equality cases in (AFI), given in terms of the support of the mixed measure $\Su(\stdb,\collc,\cdot)$, we thus obtain the following result, which establishes Schneider's conjecture \cite[Conjecture 7.6.13]{Schneider} for the class of polyoids (or smooth bodies).

\begin{theorem}
\label{corfin2}
Let $K,L\in \K^n$,  and let $\collc = (C_1, \ldots, C_{n-2})$ be a supercritical $(n-2)$-tuple of polyoids or smooth convex bodies in $\R^n$. Assume that  $\V(K,L,\collc)>0$. Then {\rm (AFI)} holds with equality if and only if there are $a>0$ and $x\in\R^n$ such that
$$h_K=h_{aL+x}\quad\text{ on } \ext (B^n,\collc).$$
\end{theorem}

In the special case where $C_1,\ldots,C_{n-2}$ are all smooth, each unit vector is  $(B^n,\collc)$-extreme and therefore $K$ and $L$ are homothetic (see \cite[Thm.7.6.8]{Schneider}).
As another consequence of Theorem \ref{thm:suppChar}, we  obtain the following partial confirmation of a conjecture on the monotonicity of mixed volumes (see \cite[Conjecture A$^{\prime}$]{Schneider1988}).

\begin{theorem}
  \makeatletter\def\@currentlabelname{Support Characterization Theorem}\makeatother\label{thm:suppChar2}
  Let $K,L\in\K^n$ satisfy $K\subseteq L$. Let  $\collc = (C_1, \ldots, C_{n-1})$ be an $(n-1)$-tuple of polyoids (or smooth convex bodies provided at least one of the bodies $C_i$ is smooth and strictly convex) in $\R^n$. Then equality holds in
  $$
  \V(K,\collc)\le \V(L,\collc)
  $$
  if and only if
  \begin{equation}\label{eq:key2}
h_K=h_L \quad\text{on }  \cextdirs
  .\end{equation}
\end{theorem}
Condition \eqref{eq:key2} is expressed by saying that $K$ and $L$ have the same $\collc$-extreme supporting hyperplanes.

In order to show relation \eqref{eq:key}, which is the main result of this paper, we prove two inclusions. Both inclusions require various preparations and involve new ideas. The main task is to prove the result in the case where   $C_1,\ldots,C_{n-1}$ are polyoids with generating measures $\mu_1,\ldots,\mu_{n-1}$. In order to show that $\csupp \subseteq \cl\cextdirs$, we express in a first step the support of the mixed area measure of $\collc$  as the closure of the union of the extreme normal vectors $\ext \mathbfcal{P}$ of all $(n-1)$-tuples $\mathbfcal{P}$ of polytopes in the support of $\mu_1\otimes\cdots\otimes\mu_{n-1}$. A main tool is Theorem \ref{thm:suppInt} which applies to more general bodies than polyoids. In a second step, we  provide in Section \ref{sec:3} information about projections of touching spaces (the linear subspaces orthogonal to the better known touching cones) and projections of polyoids and their generating measures. In Section \ref{sec:4} we develop a method to characterize what it means that the touching space of a convex body, and in particular of a polyoid, is trivial. These ingredients are combined in Section \ref{sec:7} to complete the proof of the inclusion ``$\subseteq$''. In fact, our arguments for the inclusion ``$\subseteq$'' apply to a formally larger class of convex bodies which we called macroids in \cite{HugReichert23+}, see Proposition \ref{prop:suppChar}.

For the reverse inclusion, we proceed by induction over the dimension (see Section \ref{sec:7}). A natural ingredient in the argument is a reduction formula that relates the mixed area measures of convex bodies, where some of these bodies are contained in a subspace $E$, to the mixed area measure of the remaining bodies, projected to the orthogonal subspace $E^\perp$ (see Secton \ref{sec:2}). A crucial new idea to make the induction work is to reduce the complexity of a polyoid $M$, which has a nontrivial touching space in direction $u$, by a construction we call pruning. It ultimately allows us to replace $M$ locally by a lower dimensional witness polytope $\Re(M,u)$ which can be used in place of $M$ to explore the support of a mixed area measure involving $M$. Motivating examples for the construction of such a polytope and the crucial Witness Lemma \ref{thm:pruningLemma} are contained in Section \ref{sec:5}. It is this part of the argument for the inclusion ``$\supseteq$'' which inhibits the extension of Theorem \ref{thm:suppChar} to macroids.   Another ingredient for the induction is provided in Section \ref{sec:6}. It finally allows us in the induction step to replace, for a given direction $u$, some of the polyoids  by their associated witness polytopes. The required Switching Lemma \ref{thm:critSwitching} is based on concepts of criticality that are discussed in Section \ref{sec:2} and have already  proved to be essential in recent work by Shenfeld and van Handel \cite{SvH23+}.

\section{Preparations}\label{sec:2}

We work in Euclidean space $\R^n$ with scalar product $\langle\cdot\,,\cdot\rangle$, norm $\|\cdot\|$ and Euclidean metric $d(\cdot\,,\cdot)$. Most of the time we work with nonempty compact convex subsets of $\R^n$ (convex bodies) and denote the space of all convex bodies in $\R^n$ by $\K^n$, together with the Hausdorff metric. We denote by $\Pt^n$ the subset of $\K^n$ consisting of polytopes. It is useful to consider some basic operations and concepts from convex geometry also for non-convex sets.  This is straightforward for  the Minkowski (i.e. vector) sum or Minkowski combinations with real coefficients of arbitrary subsets of $\R^n$. For $n\in \N$, we set $[n]:=\{1,\ldots,n\}$.

If $\varnothing\neq A\subseteq\R^n$, we denote by $\vspan A$ the (linear) span and by $\pspan A:=\vspan(A-A)$ the linear subspace parallel to the affine span  of $A$. Then $\dim A \coloneqq \dim \pspan A$ is the dimension of $A$. We write $\relint B$ for the relative interior of a convex set $B\subseteq\R^n$. For $x,y\in\R^n$,  the segment connecting $x$ and $y$ is denoted by $[x,y]$ (which equals the convex hull of $\{x,y\}$).

The support function of a subset $A\subseteq \R^n$ is  $ h_A \colon \R^n \to [-\infty, \infty]$, $  u \mapsto \sup\set*{ \langle x, u\rangle \given x\in A}$ and the support set of $A$ in direction $u \in \R^n \elminus0$ is
  \[
    F(A, u) \coloneqq \set*{ x \in A \given \left<x, u\right> = h_A(u) },\]
which can be the empty set. For a convex body $A$ and $u \in \R^n \elminus0$, the support set $F(A,u)$ is again a convex body.

\subsection{Faces and touching spaces}

We follow Schneider \cite[p.~16]{Schneider} in defining a face of a convex set $A\subseteq\R^n$ as a convex subset $F \subseteq A$ with the following property: If $x, y \in A$ and $F\cap\relint [x, y]\neq\varnothing$, then $[x, y] \subseteq F$.

Several useful properties of faces of nonemtpy closed convex sets are provided in \cite[Sect.2.1]{Schneider} and will be used in the following. In particular, for a polytope $P$ a set $\varnothing\neq F\subset P$  is a face of $P$ if and only if it is a support set. Note that ``$\subset$'' means strict inclusion.

Next we collect and complement some definitions from \cite[p.~85]{Schneider}. As usual, for a subset $A\subset\R^n$ we set $A^\perp:=\set*{v\in \R^n\given \langle v, a\rangle=0 \text{ for }a\in A}$ (which equals the orthogonal complement of $\vspan A$). For a vector $u\in\R^n$, we set $u^\perp:=\{u\}^\perp$.

\begin{definition}
\label{def:touching}{\rm
  Let $K$ be a convex body contained in some linear subspace $V \subseteq \R^n$.
  The set of common outer normal vectors (including $0$) of some set $S \subseteq K$ is
  \[
    N_V(K, S) \coloneqq \set*{ u \in V\elminus0 \given S \subseteq F(K, u) } \cup \set*{0} \subseteq V
  \]
  and called the \emph{normal cone of $K$ at $S$ (in $V$)}.

  If $u \in V\elminus0$, then $N_V(K, F(K, u))$ is a closed convex cone containing $u$. As such, it has a unique face $T_V(K, u)$ such that $u \in \relint T_V(K, u)$. This face is called the \emph{touching cone of $K$ in direction $u$}.

  The space $\TS_V(K, u) \coloneqq V \cap T_V(K, u)^\perp$ is called the \emph{touching space of $K$ in direction $u$}.

  In case of $V = \R^n$, we write $N \coloneqq N_{\R^n}, T \coloneqq T_{\R^n}$ and $\TS \coloneqq \TS_V$.}
\end{definition}

The following definition of extreme normal directions for an $(n-1)$-tuple of convex bodies in $\R^n$ can be easily seen to be equivalent to the
definition given in \cite[p.~87]{Schneider} by means of \cite[Lem.~5.1.9]{Schneider}, applied in $u^\perp$ for some $u\in\stdsph$.

\begin{definition}
\label{def:extremeNormal}{\rm
  If $n \ge 1$ and
  $\collc = (C_1, \ldots, C_{n-1})$ is a tuple of convex bodies in $\R^n$, then $u\in\stdsph$ is said to be a \emph{$\collc$-extreme (normal) vector} if there are  one-dimensional linear subspaces of $\TS(C_i, u)$, for $i\in[n-1]$, with linearly independent directions.
  The set of all $\collc$-extreme normal vectors is denoted by $\ext \collc$.}
  \end{definition}

  \begin{remark}{\rm
  In the situation of Definition \ref{def:extremeNormal}, $u\in \ext\collc$ if and only if
  \begin{equation}\label{eq:condiextreme}
    \dim \sum_{i \in I} \TS(C_i, u) \ge \abs{I}\quad \text{for }I \subseteq [n-1],
  \end{equation}
  where the empty sum is understood as the trivial vector space. For the equivalence of this condition with Definition \ref{def:extremeNormal}, see \cite[Thm.~5.1.8]{Schneider}.

  With the notation introduced later in Definition \ref{def:crit}, condition \eqref{eq:condiextreme} will be expressed by writing
  \[
    \V\prn*{\TS(C_1, u), \ldots, \TS(C_{n-1}, u)} > 0
  ;\]
  see also the more general Lemma \ref{thm:critIndep}.
  }
  \end{remark}

The facial structure, touching cones and touching spaces of polytopes are reasonably well-understood. In Lemmas \ref{thm:facialStability} and \ref{thm:touchingConePoly} we   provide some related information that will be needed in the sequel.

\begin{lemma}[Facial stability]\label{thm:facialStability}
  Let $P = \conv\set*{ v_1, \ldots, v_\ell } \in \Pt^n$ and $u \in \R^n \elminus0$. Consider the set
  $I_u \coloneqq \set*{ i \in [\ell] \given v_i \in F(P, u) }$.

  Then there is an $\varepsilon \in (0, \norm{u})$ such that
  for all $v, w_1, \ldots, w_\ell \in \R^n$ with $d(u, v) < \varepsilon$ and $d(v_i, w_i) < \varepsilon )$, $Q \coloneqq \conv\set*{ w_1, \ldots, w_\ell }$ satisfies $F(Q, v) \subseteq \conv\set*{ w_i \given i \in I_u }$.

  Equality holds if and only if additionally $v \in \pspan \set*{ w_i \given i \in I_u }^\perp$.
\end{lemma}
\begin{proof}
  Note that $I_u$ is not empty.

  If $I_u = [\ell]$, then the first claim follows from $F(Q, v) \subseteq Q$. Now assume that $[\ell]\setminus I_u \ne \varnothing.$
  Let $\norm{u} > \varepsilon > 0$ and $v, w_1, \ldots, w_\ell \in V$ with $d(u, v) < \varepsilon $ and $d(v_i, w_i) < \varepsilon $. Then $v \ne 0$. Define convex bodies
  \[
    Q \coloneqq \conv\set*{ w_1, \ldots, w_\ell }
  \]
  and
  \[
    P' \coloneqq \conv\set*{ v_i \given i \in [\ell] \setminus I_u }, \quad Q' \coloneqq \conv\set*{ w_i \given i \in [\ell] \setminus I_u }
  .\]
  Note that $Q$ and $Q'$ depend on $v, w_1, \ldots, w_\ell$. It holds
  \[
    h(P', u) = \max_{i \in [\ell]\setminus I_u} \left<v_i, u\right> < \max_{i \in [\ell]} \left<v_i, u\right> = h(P, u) = h(F(P, u), u).
  \]
  By continuity in $u$ and $(v_i)_{i \in [\ell]}$ of the left and right term of the inequality, we can choose $\varepsilon > 0$ such that for all $v, w_1, \ldots, w_\ell \in V$ with $d(u, v) < \varepsilon $ and $d(v_i, w_i) < \varepsilon $,
  \[
    h(Q', v) < h(Q, u)
  .\]
  So $F(Q, v) \subseteq \conv\set*{ w_i \given i \in I_u }$, remembering that $F(Q, v)$ is spanned by vertices of $Q$ \cite[Theorem~1.19]{Hug}.

  If equality holds, then $v \in \prn*{\pspan F(Q, v)}^\perp = \prn*{\pspan \set*{ w_i \given i \in I_u }}^\perp$. Conversely, assume $v \in \prn*{\pspan\set*{ w_i \given i \in I_u }}^\perp$. Since $F(Q, v)$ is spanned by vertices $w_i$ with $i\in I_u$ and nonempty, we may assume that $1 \in I_u$ and $\left<w_1, v\right> = h(Q, v)$. For $v \in \prn*{\pspan\set*{ w_i \given i \in I_u }}^\perp$, we get
  \[
    \left<w_i, v\right> = \left<w_1, v\right> +  \left<w_i - w_1, v\right>  = h(Q, v) \quad \text{for all $i \in I_u$}
  ,\]
  hence $w_i\in F(Q,v)$, and therefore also
   $ \conv\set*{ w_i \given i \in I_u }\subseteq F(Q, v) $ holds.
\end{proof}

The next lemma should be compared to \cite[(2.26)]{Schneider}.

\begin{lemma}
\label{thm:touchingConePoly}
  Let $P \in \Pt^n$ be a polytope and $u \in \R^n\elminus0$.
  Then
  $$T(P, u) = N(P, F(P, u))\quad\text{and}\quad\TS(P, u) = \pspan F(P, u).$$
\end{lemma}
\begin{proof}
If $v\in  N(P, F(P, u))$, then $F(P,u)\subseteq F(P,v)$, and thus $v\in (\pspan F(P, u))^\perp$. Hence $N(P, F(P, u))\subseteq(\pspan F(P,u))^\perp$ and therefore
  $$\pspan F(P, u) \subseteq N(P, F(P, u))^\perp \subseteq T(P, u)^\perp = \TS(P, u).$$

  By Lemma \ref{thm:facialStability}, there is an open neighborhood $U \subseteq \R^n\elminus0$ of $u$ such that for all $v \in U$,
  \[
    F(P, v) \subseteq F(P, u)
  \]
  and such that for all $v \in U' \coloneqq U \cap \prn*{\pspan F(P, u)}^\perp$, even equality holds.

  So $U' \subseteq N(P, F(P, u)) \subseteq \prn*{\pspan F(P, u)}^\perp$. But $U'$ is open in $\prn*{\pspan F(P, u)}^\perp$, so that $u \in \relint N(P, F(P, u))$. By definition of $T(P, u)$,
  \[
    T(P, u) = N(P, F(P, u))
  \]
  and hence (by the preceding argument)
  \[
    \prn*{\pspan F(P, u)}^\perp = \vspan N(P, F(P, u)) = \vspan T(P, u)
  .\]
  Thus we get $\pspan F(P, u) = T(P, u)^\perp = \TS(P, u)$.
\end{proof}

\subsection{Mixed volumes and mixed area measures}

See \cite{Schneider,Hug} for an introduction to mixed volumes and mixed area measures of convex bodies or differences of support functions of convex bodies. We start with some simple comments and conventions.

\medskip

\noindent
\textbf{Conventions concerning tuples of sets}

\medskip

\noindent
  Most of the time, the ordering of a tuple will not be relevant for our purposes. This is why a \emph{subtuple} of a tuple $\mathbfcal{A} = (A_1, \ldots, A_\ell)$, $\ell\in\N_0$, will denote any tuple $\mathbfcal{B}$ that is a prefix of a permutation of $\mathbfcal{A}$. The notation for this situation is $\mathbfcal{B} \le \mathbfcal{A}$.

  Every set $I \subseteq [\ell]$ can be uniquely written as $I = \set*{i_1, \ldots, i_m}$ such that $m\in\N_0$ and $(i_j)_{j \in [m]}$ is strictly increasing in $j$. Then we assign to $I$ a subtuple of $\mathbfcal{A}$,
  \[
    \mathbfcal{A}_I \coloneqq (A_{i_1}, \ldots, A_{i_m}) \le \mathbfcal{A}
  .\]

  The \emph{span} of a tuple of \emph{nonempty} sets $\mathbfcal{A} = (A_1, \ldots, A_\ell)$ with $A_i\subseteq\R^n$ is
  \[
    \pspan \mathbfcal{A} \coloneqq \pspan \sum_{i=1}^\ell A_\ell = \sum_{i=1}^\ell \pspan A_\ell
  ,\]
  where $\sum_{i=1}^\ell A_\ell \coloneqq  \set*{0}$ if $\ell = 0$.

  The \emph{dimension} of a tuple means the dimension of its affine span, that is, $\dim \mathbfcal{A}\coloneqq\dim\pspan \mathbfcal{A}$.
  The \emph{size} of a tuple $\mathbfcal{A}$ is the number of its components and is written as  $\abs{\mathbfcal{A}} \coloneqq \ell$.

  Whenever tuples of sets are nested into other tuples, we will omit brackets as convenient. For example, if $C, D$ are sets and $\mathbfcal{A} = (A_1, \ldots, A_\ell)$ is a tuple of sets, then
  \[
    (C, \mathbfcal{A}, D) \coloneqq (C, A_1, \ldots, A_\ell, D)
  \]
  and therefore, for example, if the right term is well-defined,
  \[
    \V(C, \mathbfcal{A}, D) = \V(C, A_1, \ldots, A_\ell, D)
  .\]

  If $\mathbfcal{A}, \mathbfcal{B}$ are tuples, we also write
  \[
    \mathbfcal{A} + \mathbfcal{B} \coloneqq (\mathbfcal{A}, \mathbfcal{B})
  ,\]
  using the nested-tuple convention as just described.

  If $k\in\N$, then for arbitrary $X$ (being a set, a measure, \ldots), $X[k]$ denotes the tuple consisting of $k$ copies of $X$, that is,
  $X[k] \coloneqq  (X, \ldots, X) $. As usual we set $$S_{n-1}(K,\cdot)\coloneqq S(K[n-1],\cdot)\quad \text{ for } K\in\K^n.
  $$

  If $f\colon A \to B$ is a function and $\mathbfcal{A} = (A_1, \ldots, A_\ell)$ is a tuple of elements or subsets of $A$, then we write
  \[
    f(\mathbfcal{A}) = f(A_1, \ldots, A_\ell) \coloneqq \prn*{f(A_1), \ldots, f(A_\ell)}
  .\]

  If $\mathbfcal{A} = (A_1, \ldots, A_\ell)$ is a tuple and $r \in [\ell]$, the tuple obtained from $\mathbfcal{A}$ by removing the $r$-th entry (i.e.\ $A_r$) is denoted by $\mathbfcal{A}_{\setminus r}$.

\begin{remark}
\label{rem:empty}{\rm
 For the discussion of mixed volumes and area measures it is usually assumed that $n\ge 1$ (or even $n\ge 2$). In view of induction arguments in the following, we set
  \[
    \V() := \V_0(\set*{0}) \coloneqq \mH{0}(\set*{0}) = 1
  ,\]
  where $\mH{0}$ is the zero-dimensional Hausdorff measure (counting measure).
Moreover, for $n=1$ we define $\Su()$ as the counting measure on $S^0 = \set*{-1, 1}$. Then e.g. relation \eqref{eqex} remains true.

  These definitions are consistent with the inductive definitions of volume and surface in \cite[Definition~3.2]{Hug}.}
\end{remark}

\begin{remark}
{\rm
  In order to simplify notation, we use the following conventions.
\begin{enumerate}[{(1)}]
\item   Let $\mu(\mathbfcal{C})$ be a measure which depends on a parameter $\mathbfcal{C}$. Then we write $\mu(\mathbfcal{C}, \cdot)$ or $\mu_{\mathbfcal{C}}(\cdot)$ as shorthands for $\mu(\mathbfcal{C})(\cdot)$.
\item  Sometimes it is useful to pass the support function $h_K$ instead of the convex body $K \in\K^n$ to $\Su$ or $\V$, i.e., write $\V(h_K, \mathbfcal{C})$ instead of $\V(K, \mathbfcal{C})$. Using this convention, $\V$ (and $\Su$) can be extended to multilinear functions taking $n$ (or $n-1$) differences of support functions. For example,
  \[
    \V(h_K - h_L, \mathbfcal{C}) \coloneqq  \V(K, \mathbfcal{C}) - \V(L, \mathbfcal{C})
  .\]
  \item In the following we write $\V$ for the mixed volume in $\R^n$, but we use the same symbol for the mixed volume in a subspace (the number of arguments already provides the relevant information). By the translation invariance of mixed volumes, the mixed volume of convex bodies lying in parallel subspaces is well-defined.
  \end{enumerate}}
\end{remark}

The mixed area measure of an $(n-1)$-tuple of polytopes can be written as a finite sum of weighted Dirac measures and the point mass (weight) of each atom is given as a mixed volume. We recall this relation in the remark below since it will be used in the following and a related (more general) result for general convex bodies is stated as Lemma \ref{thm:areaVol}.

\begin{remark}
\label{rem:maPoly}{\rm
  Let $P_1, \ldots, P_{n-1}\in \Pt^n$ and $P \coloneqq P_1 + \cdots + P_{n-1}$. Then the mixed area measure of $P_1, \ldots, P_{n-1}$ is a weighted sum of Dirac measures, that is,
  \[
    \Su(P_1, \ldots, P_{n-1},\cdot) = \sum_{u \in \mathcal{N}_{n-1}(P)} \V(F(P_1, u), \ldots, F(P_{n-1}, u)) \delta_u
  ,\]
  where $\mathcal{N}_{n-1}(P)$ is the set of all $u\in \mathbb{S}^{n-1}$ with $\dim F(P,u)=n-1$
  (see \cite[(4.2)]{Hug}).}
\end{remark}

We will end this discussion by recalling a useful result which relates the mixed area measure $S_{n-1}(K;\cdot)$ of the $(n-1)$-tuple $(K,\ldots,K)$ to the (localized) $(n-1)$-dimensional Hausdorff measure $\mathcal{H}^{n-1}$  of the topological boundary $\partial K$  of an $n$-dimensional convex body $K$ in $\R^n$.

\begin{definition}
{\rm
  Let $n\ge 1$, $K \in \K^n$ a convex body and $\omega \subseteq \stdsph$ a set. Then
  \[
    \tau(K, \omega) \coloneqq \bigcup_{u \in \omega} F(K, u)
  \]
  is called the \emph{reverse spherical image of $K$ at $\omega$ } (compare \cite[p.~88]{Schneider}).}
\end{definition}

\begin{lemma}\label{thm:areaHaus}
  Let $n \ge 1$. For every $n$-dimensional convex body $K \in \K^n$ and every Borel measurable set $\omega \subseteq \stdsph$,
  \[
   \Su_{n-1}(K,\omega)  = \mH{n-1}(\tau(K, \omega))
  .\]
\end{lemma}
\begin{proof}
  See \cite[Theorem~4.2.3]{Schneider} or \cite[Thm.~4.8]{Hug}.
\end{proof}

Lemma \ref{thm:areaHaus}  in combination with the well-known (diagonality) Lemma \ref{thm:multilin} has many applications, such as Lemmas \ref{thm:maTau} and \ref{thm:areaVol}.

\begin{lemma}
  \makeatletter\def\@currentlabelname{Diagonality Lemma}\makeatother\label{thm:multilin}
  Let $f, g\colon (\K^n)^k \to \R$ be functionals that are symmetric and multilinear (i.e.\ Minkowski additive and positively homogeneous in each of their $k \in \N_0$ components) and let $\collc = (C_1, \ldots, C_k)$ be a tuple of convex bodies in $\R^n$. If for all choices of $\lambda = (\lambda_1,\ldots,\lambda_k)\in [0, \infty)^k$ the convex body
  \[
    L_\lambda \coloneqq \sum_{i=1}^k \lambda_i C_i
  \]
  satisfies $f(L_\lambda[k]) = g(L_\lambda[k])$, then $f(\collc) = g(\collc)$.
\end{lemma}

The following lemma states that the mixed area measures are locally determined, which will be crucial for the proof of Lemma \ref{thm:realPruningLemma} (and it will be used in the discussion of some of the examples). For the area measures of a single convex body (and Euclidean balls), the corresponding simple fact is well known (see \cite[Note 11 for Sect.~4.2]{Schneider}).

\begin{lemma}
\label{thm:maTau}
  Let $n \ge 1$. Let $\collc = (C_1, \ldots, C_{n-1}), \mathbfcal{D} = (D_1, \ldots, D_{n-1})$ be tuples of convex bodies in $\R^n$, and let  $\omega \subseteq \stdsph$ be a Borel set such that
  \[
    \tau(C_i, \omega) = \tau(D_i, \omega),\quad i \in [n-1].
  \]
  Then
  \[
    \Su(\collc)(\omega) = \Su(\mathbfcal{D})(\omega)
  .\]
\end{lemma}
\begin{proof}
  The case $n = 1$ follows from the fact that $\collc = \mathbfcal{D}$ are empty tuples. We will prove the theorem for the case that $n \ge 2$ and $C_i = D_i$ for $i \ne 1$. This allows one to replace $C_1$ by $D_1$, yielding
  \begin{align}\label{eq:prepequal}
    \Su(C_1, C_2, \ldots, C_{n-1})(\omega) = \Su(D_1, C_2, \ldots, C_{n-1})(\omega)
  .\end{align}
  Using symmetry of $\Su$, we can afterwards replace $C_2$ by $D_2$,  and so on until we have replaced all $C_i$ by $D_i$.

We start with a preparatory remark.  Let $K \in \K^n$, $\omega\subseteq\mathbb{S}^{n-1}$ and $u\in \omega$. We show that $F(K,u)=F(\tau(K, \omega), u)$. First, observe that $F(K, u) \subseteq \tau(K, \omega) \subseteq K$. Hence, $h_{\tau(K, \omega)}(u) = h_K(u)$ and
  \begin{align*}
    F(K, u) &= \set*{ x \in K \given \left<x, u\right> = h_K(u) } = \set*{ x \in \tau(K, \omega) \given \left<x, u\right> = h_{\tau(K, \omega)}(u) }\\
    & = F(\tau(K, \omega), u)
  ,\end{align*}
  where we again used that $F(K, u) \subseteq \tau(K, \omega)$.

  By Minkowski additivity of the mixed area measure in its first component, it suffices to show that \eqref{eq:prepequal} holds when $C_1, D_1$ are full-dimensional. To see this, replace $C_1$ by $C_1 + \stdb$ and $D_1$ by $D_1 + \stdb$ and note that $\tau(C_1 + \stdb,\omega)=\tau(D_1 + \stdb,\omega)$, since  by the preparatory remark for any $u\in \omega$ we have
  \begin{align*}
  F(C_1+B^n,u)&=F(C_1,u)+F(B^n,u)=F(\tau(C_1,\omega),u)+F(B^n,u)\\
  &=F(\tau(D_1,\omega),u)+F(B^n,u) =
  F(D_1,u)+F(B^n,u)=F(D_1+B^n,u).
  \end{align*}

  For every $(\lambda_i)_{i\in[n-1]}\in[0, \infty)^{n-1}$, we claim that
  \begin{align}\label{eq:mt1}
    \Su_{n-1}\prn*{\sum_{i=1}^{n-1}\lambda_i C_i}(\omega) = \Su_{n-1}\prn*{\lambda_1D_1 + \sum_{i=2}^{n-1}\lambda_i C_i}(\omega)
  .\end{align}
  If this holds, Lemma \ref{thm:multilin} will show
  \[
    \Su(C_1, C_2, \ldots, C_{n-1})(\omega) = \Su(D_1, C_2, \ldots, C_{n-1})(\omega)
  .\]

  If $\lambda_1 = 0$,  \cref{eq:mt1} clearly holds. Otherwise, $\sum_{i=1}^{n-1}\lambda_i C_i$ and $\lambda_1 D_1 + \sum_{i=2}^{n-1}\lambda_i C_i$ are full-dimensional and  by Lemma \ref{thm:areaHaus} and the definition of $\tau$ it suffices to show that, for all $u \in \omega$,
  \begin{align*}
    F\prn*{\sum_{i=1}^{n-1}\lambda_i C_i, u} &= \sum_{i=1}^{n-1}\lambda_i F(C_i, u) \overset{(!)}{=} \lambda_1 F(D_1, u) + \sum_{i=2}^{n-1}\lambda_i F(C_i, u) \\
    &= F\prn*{\lambda_1 D_1 + \sum_{i=2}^{n-1}\lambda_i C_i, u}
  ,\end{align*}
  where we used at (!) that by the preparatory  remark  and the assumption we have
  \[
    F(C_1, u) = F(\tau(C_1, \omega), u) = F(\tau(D_1, \omega), u) = F(D_1, u)
  ,\]
  concluding the proof.
\end{proof}

The next lemma is a simple consequence of Lemma \ref{thm:maTau}, but we will not need it in the current work.

\begin{lemma}
\label{thm:areaVol}
  Assume $n \ge 1$. Let $K_1, \ldots, K_{n-1}\subset\R^n$ be convex bodies and $u \in \stdsph$. Then
  \[
    \Su(K_1, \ldots, K_{n-1})(\set*{u}) = \V(F(K_1, u), \ldots, F(K_{n-1}, u))
  .\]
\end{lemma}
\begin{proof}
  By multilinearity and symmetry of $\Su$ and $\V$ and linearity of $F$, it suffices by Lemma \ref{thm:multilin} to prove the statement for $K_1 = \cdots = K_{n-1}$, i.e. to prove that
  \[
    \Su_{n-1}(K_1)(\set*{u}) = \V_{n-1}(F(K_1, u))
  ,\]
  where $\V_{n-1}$ is the volume (intrinsic volume of order $n-1$) in an $(n-1)$-dimensional subspace of $\R^n$.
  Consider the truncated convex cone
  \[
    C \coloneqq \set*{ x \in \stdb \given  \left<x, u\right> \le -\frac{1}{2}\norm{x} }
  ,\]
  which is a full-dimensional convex body satisfying $F(C, u) = \set*{0}$. So
  \[
    \tau(K_1, \set*{u}) = F(K_1, u) = F(C + K_1, u) = \tau(C + K_1, \set*{u})
  .\]
  By Lemmas \ref{thm:maTau} and \ref{thm:areaHaus} and since $\dim(C+K_1)=n$, it follows that
  \[
    \Su_{n-1}(K_1)(\set*{u}) = \Su_{n-1}(C + K_1)(\set*{u}) = \mH{n-1}(F(C + K_1, u)) = \V_{n-1}(F(K_1, u))
  ,\]
  which completes the argument.
\end{proof}

\subsection{Reduction formulas}

We will use dimensional induction to prove assertions about mixed area measures. To succeed in this endeavor, we have to relate mixed area measures in $\R^n$ to mixed area measures in subspaces. By using basic integral geometry,  the following two reduction formulas can be obtained.

Recall that we write $\V$ for the mixed volume in $\R^n$ and use the same symbol for the mixed volume in a subspace (the number of arguments already provides the relevant information). For a linear subspace $L\subseteq\R^n$, the orthogonal projection to $L$ is denoted by $\pi_L:\R^n\to L$.

\begin{lemma}
  \makeatletter\def\@currentlabelname{Reduction Formula}\makeatother\label{thm:mvRed}
  Let $\collc = (C_1, \ldots, C_{n})$ be a tuple of convex bodies in $\R^n$, and let $k \in [n]\cup\set*{0}$ be such that $\pspan \collc_{[k]}$ is contained in a linear subspace $E \subseteq \R^n$ of dimension $k$.
  Then
  \[
    \begin{pmatrix}
      n \\ k
    \end{pmatrix}
    \V(\collc) = \V(\collc_{[k]}) \cdot \V(\pi_{E^\perp}(\collc_{[n]\setminus[k]}))
  .\]
\end{lemma}
\begin{proof}
  The cases $k \in \set*{0, n}$ are trivial. For the remaining cases,
  use the translation invariance of $\V$ and apply \cite[Theorem~5.3.1]{Schneider}.
\end{proof}

In dealing with mixed area measures, we will indicate by our notation in which subspace the measure is applied. For an $\ell$-dimensional linear subspace $L\subset\R^n$, $\ell\ge 1$, we write $\Su_L$ for the mixed area measure in $L$, which is evaluated at $\ell-1$ convex bodies in $L$ and Borel subsets of $\mathbb{S}^{n-1}\cap L$. Moreover, we define $\Su_L'$ as the Borel measure on $\mathbb{S}^{n-1}$ defined by $$\Su_L'(C_1,\ldots,C_{\ell-1})(\omega)\coloneqq S_L(C_1,\ldots,C_{\ell-1})(\omega\cap L)$$ for convex bodies $C_1,\ldots,C_{\ell-1}\subset L$ and Borel sets $\omega\subseteq \mathbb{S}^{n-1}$.

The following proposition will be essential for the proof of our main result in Section \ref{sec:7}.

\begin{proposition}
  \makeatletter\def\@currentlabelname{Reduction Formula}\makeatother\label{thm:maRed}\label{thm:splittingSc}
  Assume  $n \in\N$. Let $\collc = (C_1, \ldots, C_{n-1})$ be a tuple of convex bodies in $\R^n$, and let $k \in [n-1]\cup\set*{0}$ be such that $\pspan \collc_{[k]}$ is contained in a linear subspace $E \subseteq \R^n$ of dimension $k$.
  Then
  \[
    \begin{pmatrix}
      n-1 \\ k
    \end{pmatrix}
    \Su(\collc) = \V(\collc_{[k]}) \cdot \Su'_{E^\perp}(\pi_{E^\perp}(\collc_{[n-1]\setminus[k]}))
  .\]
  In particular, if $\dim \collc_{[k]} < k$, then $\Su(\collc) = 0$.
\end{proposition}
\begin{proof}
  The case $k=0$ is trivial. The assertion for $k=n-1$ is clear for polytopes (see Remark \ref{rem:maPoly}), the general case follows by approximation. So we can assume that $n \ge 3$ and $k \in [n-2]$.
  Let $C_n \in \K^n$. Then by Lemma \ref{thm:mvRed},
  \[
    \begin{pmatrix}
      n \\ k
    \end{pmatrix}
    \V(C_1, \ldots, C_n) = \V(C_1, \ldots, C_k) \cdot \V(\pi_{E^\perp}(C_{k+1}, \ldots, C_n))
  .\]
  Expressing the mixed volumes by mixed area measures, we obtain
  \begin{align*}
   & \begin{pmatrix}
      n \\ k
    \end{pmatrix}
    \frac{n - k}{n}\int h_{C_n} \diff \Su(C_1, \ldots, C_{n-1}) \\
    &= \V(C_1, \ldots, C_k) \cdot \int h_{\pi_{E^\perp} C_n} \diff\Su_{E^\perp}(\pi_{E^\perp}(C_{k+1}, \ldots, C_{n-1})
  .\end{align*}

  Noting that $h_{\pi_{E^\perp} C_n} = h_{C_n}$ on $E^\perp$, we find that
  \begin{align*}
    \int h_{\pi_{E^\perp} C_n} \diff\Su_{E^\perp}(\pi_{E^\perp}(C_{k+1}, \ldots, C_{n-1})) &= \int h_{C_n} \diff\Su_{E^\perp}(\pi_{E^\perp}(C_{k+1}, \ldots, C_{n-1})) \\
      &= \int h_{C_n} \diff\Su_{E^\perp}'(\pi_{E^\perp}(C_{k+1}, \ldots, C_{n-1}))
  \end{align*}
  and conclude that
  \begin{align*}
   & \begin{pmatrix}
      n - 1 \\ k
    \end{pmatrix}
    \int h_{C_n} \diff\Su(C_1, \ldots, C_{n-1}) \\
    &= \V(C_1, \ldots, C_k) \cdot \int h_{C_n} \diff\Su'_{E^\perp}(\pi_{E^\perp}(C_{k+1}, \ldots, C_{n-1})
  .\end{align*}
Because $C_n$ is an arbitrary convex body and differences of support functions are dense in $C(\stdsph)$, the claim follows.
\end{proof}

\subsection{Criticality}

Criticality is a useful concept that describes dimensionality conditions on arrangements of convex bodies.
Shenfeld and van Handel \cite{SvH23+} employed criticality in their investigation of equality cases in the Alexandrov--Fenchel inequality for polytopes. We will slightly deviate from their terminology in that we call ``semicritical'' what they called ``subcritical'', and we say ``subcritical'' to describe a situation which is ``not critical''.

The most elementary occurrence and motivation for the terminology is the following result.

\begin{lemma}
\label{thm:mvVanish}
  Let $\mathbfcal{C} = (K_1, \ldots, K_n)$ be a tuple of convex bodies in $\R^n$. Then the following are equivalent:
  \begin{enumerate}[{\rm (a)}]
    \item $\V(\mathbfcal{C}) > 0$.
    \item There are segments $S_i \subseteq K_i$ ($i \in [n]$) with linearly independent directions.
    \item Whenever $\mathbfcal{D} \le \mathbfcal{C}$, then $\dim \pspan \mathbfcal{D} \ge \abs{\mathbfcal{D}}$.
  \end{enumerate}
\end{lemma}
\begin{proof}
  See \cite[Theorem~5.1.8]{Schneider}.
\end{proof}

Condition (c) in Lemma \ref{thm:mvVanish} suggests the definition of a \enquote{semicritical} tuple of convex bodies. Let us recall concepts of criticality and describe some consequences.

\begin{definition}
\label{def:crit}
{\rm
  Let $\ell \in \N_0$. Let $\mathbfcal{A} = (A_1, \ldots, A_\ell)$ be a tuple of nonempty subsets of \ $\R^n$.
  Then $\mathbfcal{A}$ is called
  \begin{enumerate}[{\rm (i)}]
    \item \emph{semicritical} if for all $() \ne\mathbfcal{B} \le \mathbfcal{A}$ we have $\dim \pspan \mathbfcal{B} \ge \abs{\mathbfcal{B}}$,
    \item \emph{critical} if for all $() \ne \mathbfcal{B} \le \mathbfcal{A}$ we have $\dim\pspan\mathbfcal{B} \ge \abs{\mathbfcal{B}} + 1$,
    \item \emph{supercritical} if for all $() \ne \mathbfcal{B} \le \mathbfcal{A}$ we have $\dim\pspan\mathbfcal{B} \ge \abs{\mathbfcal{B}} + 2$,
    \item \emph{subcritical} if it is not critical.
  \end{enumerate}
  Abusing notation, we write $\V(\mathbfcal{A}) > 0$ to say that $\mathbfcal{A}$ is semicritical.}
\end{definition}

The following lemma is provided in \cite[Lem.~3.2]{HugReichert23+} (see also the preceding remarks there).

\begin{lemma}
\label{thm:critSimple}
  Let $\ell \in \N_0$,  and let $\mathbfcal{A} = (A_1, \ldots, A_\ell)$ be a tuple of  nonempty subsets of  \  $\R^n$.
  \begin{enumerate}[{\rm (1)}]
    \item Subtuples of (super-, semi-)critical tuples are also (super-, semi-)critical.\label{it:critSimple3}
    \item Supercriticality implies criticality, which implies semicriticality.\label{it:critSimple4}
    \item The empty tuple is supercritical.\label{it:critSimple6}
    \item (Super-, Semi-)Criticality is invariant under permutations of $\mathbfcal{A}$.\label{it:critSimple7}
    \item (Super-, Semi-)Criticality is invariant under simultaneous affine isomorphisms and argumentwise translations.\label{it:critSimple8}
    \item (Super-, Semi-)Criticality is preserved if the sets in $\mathbfcal{A}$ are replaced by supersets.\label{it:critSimple9}
    \item Let $\mathbfcal{A}$ be critical and $A_{\ell + 1} \subseteq \R^n$ be nonempty. Then $(A_1, \ldots, A_{\ell + 1})$ is semicritical if and only if $A_{\ell + 1}$ is at least one-dimensional.\label{it:critSimple1}
    \item Let $\mathbfcal{A}$ be supercritical and $A_{\ell + 1} \subseteq \R^n$ be nonempty. Then $(A_1, \ldots, A_{\ell+1})$ is critical if and only if $A_{\ell+1}$ is at least two-dimensional.\label{it:critSimple2}
    \item If all sets $A_i$ are full-dimensional, then $\mathbfcal{A}$ is supercritical if and only if $\ell \le n-2$ or $\mathbfcal{A} = ()$.\label{it:critSimple5}
  \end{enumerate}
\end{lemma}

The  notation `$\V(\mathbfcal{A}) > 0$' suggests that semicriticality might abide laws similar to the ones applying to mixed volumes. In particular, we might hope for some kind of reduction theorem in analogy to \nref{thm:mvRed}. As the next result shows, this hope is not in vain.

The following Lemmas \ref{thm:critRed} and \ref{thm:critAdd} will be crucial for the arguments in Sections \ref{sec:6} and \ref{sec:7}. Lemma \ref{thm:critIndep} is used in the proof of Lemma \ref{thm:critAdd}.

\begin{lemma}[Semicritical reduction]\label{thm:critRed}
  Let $\ell\in\N_0$. Let $\mathbfcal{A} = (A_1, \ldots, A_\ell)$ be a tuple of nonempty subsets of \ $\R^n$ and let $\pspan \mathbfcal{A}_{[k]}$ be contained in a linear subspace $E$ of dimension $k \in \N_0$. Then the following are equivalent:
  \begin{enumerate}[{\rm (a)}]
    \item $\V(\mathbfcal{A}) > 0$;
    \item $\V(\mathbfcal{A}_{[k]}) > 0$ and $\V(\pi_{E^\perp}(\mathbfcal{A}_{[\ell]\setminus[k]})) > 0$.
  \end{enumerate}
\end{lemma}
\begin{proof}
  After applying suitable translations, we may assume that all sets contain $0$.

  \noindent
  \enquote{$\implies$}: Assume that $\V(\mathbfcal{A}) > 0$. Then by Lemma \ref{thm:critSimple}, $\V(\mathbfcal{A}_{[k]}) > 0$.

  It remains to show the second claim. For this, let $I \subseteq [\ell]\setminus[k]$. Then using the dimension formula from linear algebra and semicriticality of $\mathbfcal{A}$,
  \begin{align*}
    \dim \pspan \pi_{E^\perp}(\mathbfcal{A})_I &= \dim \pi_{E^\perp}\prn*{\pspan \mathbfcal{A}_{I \cup [k]}}
                                      \ge \dim \pspan\mathbfcal{A}_{I \cup [k]} - \dim\ker\pi_{E^\perp}\nonumber \\&\ge \abs{I} + k - k.
  \end{align*}
\noindent
  \enquote{$\impliedby$}: Now assume that $\V(\mathbfcal{A}_{[k]}) > 0$ and $\V(\pi_{E^\perp}(\mathbfcal{A}_{[\ell]\setminus[k]})) > 0$.

  Let $I \subseteq [\ell]$ and consider the linear map
  $\Phi \colon \pspan \mathbfcal{A}_I \to \R^n$, $x \mapsto \pi_{E^\perp}(x)$.
  It satisfies
  \[
    \ker \Phi = E \cap \pspan \mathbfcal{A}_I \supseteq \pspan \mathbfcal{A}_{I \cap [k]}
  \]
  and
  \[
    \im \Phi = \pspan \pi_{E^\perp}(\mathbfcal{A}_I) = \pspan \pi_{E^\perp}(\mathbfcal{A}_{I \setminus [k]})
  .\]
  The dimension formula together with the assumption shows
  \begin{align*}
   \dim\pspan \mathbfcal{A}_I & = \dim \ker\Phi + \dim \im \Phi \\
    &\ge \dim\pspan \mathbfcal{A}_{I \cap [k]} + \dim \pspan \pi_{E^\perp}(\mathbfcal{A}_{I \setminus [k]}) \\
    &= \dim\pspan \mathbfcal{A}_{I \cap [k]} + \dim\pspan \prn*{\pi_{E^\perp}(\mathbfcal{A})}_{I \setminus [k]} \\
    &\ge \abs{I \cap [k]} + \abs{I \setminus [k]} = \abs{I}
  ,\end{align*}
  which shows that $\mathbfcal{A}$ is semicritical.
\end{proof}

Having proved the reduction Lemma \ref{thm:critRed}, we can inductively prove an analogue of Lemma \ref{thm:mvVanish}.

\begin{lemma}
\label{thm:critIndep}
  Let $\mathbfcal{A} = (A_1, \ldots, A_\ell)$ be a tuple of nonempty subsets of $\R^n$. Then the following are equivalent:
  \begin{enumerate}[{\rm (a)}]
    \item $\V(\mathbfcal{A}) > 0$.
    \item There are pairs of points $(x_i, y_i) \in A_i \times A_i$ ($i \in [\ell]$) such that the tuple $(y_i - x_i)_{i \in [\ell]}$ consists of linearly independent vectors.
  \end{enumerate}
\end{lemma}
\begin{proof}
  \enquote{$\impliedby$}: Clearly, whenever $I \subseteq [\ell]$,
  \[
    \dim \pspan \mathbfcal{A}_I \ge \dim \vspan \set*{ y_i - x_i \given i \in I } = \abs{I}
  .\]
\noindent
  \enquote{$\implies$}: We may assume that every set $A_i$ contains $0$.

  We proceed by induction over the dimension $n$. Assume that the claim is true for all dimensions smaller than $n$. Then we distinguish three cases:
  \begin{itemize}
    \item If $n = 0$, we have nothing to show since the empty family is clearly linearly independent.
    \item If $n > 0$ and the tuple is critical, let $E$ be an arbitrary $(n-1)$-dimensional linear subspace. Then $\pi_E(\mathbfcal{A})$ is still semicritical because the kernel of the projection is one-dimensional. The inductive hypothesis guarantees the existence of pairs of points $(x_i, y_i) \in A_i \times A_i$ ($i \in [\ell]$) such that $\pi_E(y_i - x_i) \in E$ are linearly independent. But then $(y_i - x_i)$ are linearly independent, too.
    \item If $n > 0$ and the tuple is subcritical, we find $\varnothing \ne I \subseteq [\ell]$ with $\dim \mathbfcal{A}_I = \abs{I}$. If $\ell=n$ and $\dim \pspan \mathbfcal{A}_I=|I|$ for all $\varnothing \ne I \subseteq [n]$, then clearly there exist points $(x_i, y_i) \in A_i \times A_i$ ($i \in [\ell]$) such that the family $(y_i - x_i)_{i \in [\ell]}$ is linearly independent. Otherwise,
    without loss of generality, $\mathbfcal{A}_I$ is a prefix of $\mathbfcal{A}$, so that $I = [k]$ for some $0 < k < n$.

      After defining the linear subspace $E := \pspan \mathbfcal{A}_{[k]}$ of dimension $k$, we can apply Lemma \ref{thm:critRed} to deduce that
      \[
        \V(\mathbfcal{A}_{[k]}), \V(\pi_{E^\perp} (\mathbfcal{A}_{[\ell]\setminus[k]})) > 0
      .\]

      Because $0 < k < n$, two applications of the inductive hypothesis yield pairs of points $(x_i, y_i) \in A_i \times A_i$ ($i \in [\ell]$) such that
      \begin{itemize}
        \item[$\bullet$] $y_1 - x_1, \ldots, y_k - x_k$ are linearly independent and
        \item[$\bullet$] $\pi_{E^\perp} (y_{k+1} - x_{k+1}), \ldots, \pi_{E^\perp}(y_\ell - x_\ell)$ are linearly independent.
      \end{itemize}
     Hence it follows that $y_1 - x_1, \ldots, y_\ell - x_\ell$ are linearly independent.
  \end{itemize}
  Since these cases are exhaustive, the proof is complete.
\end{proof}

In analogy to the additivity of the mixed volume, we obtain the following result.

\begin{lemma}[Semicritical additivity]\label{thm:critAdd}
  Let $\mathbfcal{A} = (A_1, A_2, \ldots, A_\ell)$ be a tuple of nonempty subsets of $\R^n$ and $\ell \ge 1$. Furthermore, let $A_1 = B + C$. Then the following are equivalent.
  \begin{enumerate}[{\rm (a)}]
    \item $\V(A_1, A_2, \ldots, A_\ell) > 0$.
    \item $\V(B, A_2, \ldots, A_\ell) > 0$ or $\V(C, A_2, \ldots, A_\ell) > 0$.
  \end{enumerate}
\end{lemma}
\begin{proof}
  \enquote{$\impliedby$} follows from $\pspan B, \pspan C \subseteq \pspan A_1$.

  \medskip

\noindent
  \enquote{$\implies$}: In view of Lemma \ref{thm:critIndep}, we find pairs of points $(x_i, y_i)\in A_i\times A_i$ for $i \in [\ell]$ such that the differences $y_i - x_i$ are linearly independent. In particular, $y_1 - x_1$ is not contained in $E \coloneqq \vspan \set*{y_i - x_i \given i \in [\ell]\elminus 1}$.

  We can find $b, b' \in B$ and $c, c' \in C$ such that $x_1 = b + c$ and $y_1 = b' + c'$. Then either $b' - b$ or $c' - c$ is not contained in $E$ --- we may assume that $b' - b \notin E$. But then $(b'-b,y_2-x_2,\ldots,y_{\ell}-x_{\ell})$ are linearly independent, which
yields $\V(B, A_2, \ldots, A_\ell) > 0$   via Lemma \ref{thm:critIndep}.
\end{proof}

\subsection{Support of mixed area measures}

The support of mixed area measures is the central topic of this work. This section provides some of its properties that will be needed.
In the special case of polytopes, Theorem \ref{thm:suppChar} is known and easy to verify. For the sake of completeness and to familiarize the reader with our notation, we include the argument.

\begin{lemma}\label{thm:suppCharPoly}
 Let $n\ge 1$. Let $\mathbfcal{P} = (P_1, \ldots, P_{n-1})$ be a tuple of polytopes in $\R^n$. Then
  \[
    \supp\Su(\collp) = \cl\ext\collp
  .\]
\end{lemma}
\begin{proof}
For $n=1$ the assertion is clear by our definitions. Let $n\ge 2$. By Remark \ref{rem:maPoly},
  \begin{align}\label{eq:su1}
    \Su(\collp) = \sum_{u \in \mathcal{N}_{n-1}(P_1 + \cdots + P_{n-1})} \V(F(P_1, u), \ldots, F(P_{n-1}, u)) \delta_u,
  \end{align}
  and for all $u \in \stdsph$, Lemmas \ref{thm:touchingConePoly} and \ref{thm:mvVanish} show the equivalence
  \begin{align}\label{eq:su2}
    \V(F(P_1, u), \ldots, F(P_{n-1}, u)) > 0 \, \iff \, \V(\TS(P_1, u), \ldots, \TS(P_{n-1}, u)) > 0,
  \end{align}
  the second statement by definition being equivalent to $u \in \ext\collp$.

  So if $u \in \supp\Su(\collp)$, then $\V(\TS(P_1, u), \ldots, \TS(P_{n-1}, u)) > 0$, i.e.\ $u \in \ext\collp$. Therefore, $\supp\Su(\collp) \subseteq \ext\collp \subseteq \cl\ext\collp$.

  Conversely, assume $u \in \ext\collp$. Then $\V(F(P_1, u), \ldots, F(P_{n-1}, u)) > 0$ follows from  \eqref{eq:su2}. In particular,
  \[
    \dim F\prn*{\sum_{i=1}^{n-1} P_i, u} = \dim \sum_{i=1}^{n-1} F(P_i, u) \ge n-1
  .\]
  So $u \in \mathcal{N}_{n-1}\prn*{\sum_{i=1}^{n-1} P_i}$ and
  \[
    \Su(\collp)(\set*{u}) \ge \V(F(P_1, u), \ldots, F(P_{n-1}, u)) > 0
  ,\]
  which shows that $u \in \supp\Su(\collp)$, hence  $\ext\collp \subseteq \supp\Su(\collp)$. The claim
  follows,  since $\supp\Su(\collp)$ is closed.
\end{proof}

Next we describe the support of a convex body which is defined as an integral average in terms of its support function.

\begin{theorem}
\label{thm:suppInt}
  Assume that $n \ge 2$.
  Let $C_1 \in \K^n$ be a convex body, $\collc = (C_2, \ldots, C_{n-1})$ an $(n-2)$-tuple of convex bodies and $\mu$ a finite Borel measure on $\K^n$ with bounded  support such that
  \[
    h_{C_1}(x) = \int h_K(x) \,\mu(\diff K),\quad x\in\R^n
  .\]
  Then
  $$
  \Su(C_1, \collc)=\int S(K,\collc)\, \mu(\diff K)
  $$
  and
  \[
    \supp \Su_{C_1, \collc} = \cl \bigcup_{K \in \supp\mu} \supp \Su_{K, \collc}
  .\]
\end{theorem}
\begin{proof}
  Let $A \subseteq \stdsph$ be closed. Let $d(u, A)$ denote the Euclidean distance of $u\in\mathbb{S}^{n-1}$ from $A$.
  Then the continuous function
  \[
    f_A \colon \stdsph \to [0, \infty), \quad u \mapsto d(u, A),
  \]
  satisfies $f_A^{-1}(\set*{0}) = A$.

  If $f$ is a difference of support functions, we can apply Fubini's theorem and the compactness of the support of $\mu$ to obtain
  \begin{align*}
    \int f \diff \Su_{C_1, \collc} &= \int h_{C_1} \diff \Su_{f, \collc} \\
      &= \int \int h_K(x) \, \mu(\diff K) \Su_{f, \collc}(\diff x) \\
      &= \int \int h_K(x) \,\Su_{f, \collc}(\diff x) \,\mu(\diff K) \\
      &= \int \int f(x) \,\Su_{K, \collc}(\diff x)\, \mu(\diff K).
  \end{align*}
  The same equality holds for all continuous functions $f \colon \stdsph \to \R$ by approximation, and in particular  for $f_A$ as defined above. Thus we have verified the first assertion. Now we turn to the second claim.

  \medskip

\noindent
  \enquote{$\subseteq$}:
  Set $f \coloneqq f_{\cl \bigcup_{K \in \supp\mu} \supp \Su_{K, \collc}}$. Then
  \[
    \int f \diff \Su_{C_1, \collc} = \int\int f(x) \,\Su_{K, \collc}(\diff x)\, \mu(\diff K) = 0
  .\]
  So $\Su_{C_1, \collc}(f^{-1}((0, \infty))) = 0$, concluding this direction.

  \medskip

\noindent
  \enquote{$\supseteq$}: Let $x \notin \supp \Su_{C_1, \collc}$. Because $\supp\Su_{C_1, \collc}$ is closed, it suffices to prove that $x \notin \supp\Su_{K, \collc}$ for all $K \in\supp\mu$.

  There is an open set $U \subseteq \stdsph$ with $x \in U$ such that $\Su_{C_1, \collc}(U) = 0$. Define $f \coloneqq f_{U^\complement}$.
  Then
  \[
    0 = \int f \diff \Su_{C_1, \collc} = \int\int f(z) \,\Su_{K, \collc}(\diff z)\, \mu(\diff K)
  .\]
  The integrand $\varphi \colon K \mapsto \int f(z)\, \Su_{K, \collc}(\diff z)$ is nonnegative and continuous by the continuity of $f$ and the weak continuity of the mixed area measure.
  Therefore, $\varphi(K)=0$  for $K\in\supp\mu$.
  In other words, if $K \in \supp\mu$, then $\int f \diff \Su_{K, \collc} = 0$. The integrand being nonnegative and continuous, $f$ vanishes on $\supp\Su_{K, \collc}$. Therefore,
  \[
    x \in U \subseteq (\supp\Su_{K, \collc})^\complement
  ,\]
  which was to be shown.
\end{proof}

The preceding theorem can in particular be applied in the case where $C_1$ is a polyoid, as follows from \cite[Cor.~2.9]{HugReichert23+}.

Finally, we mention a general result which states that the support of the weak limit of a sequence of measures is covered (up to taking the closure) by the supports of these measures. The proof is a straightforward consequence of the definition of weak convergence of measures.

\begin{lemma}[Support and weak convergence]\label{thm:suppEqOfTendsto}
  Let $\mu_\ell \to \mu$ be a weakly convergent sequence of finite Borel measures on a second-countable metric space $E$. Then
  \[
   \supp \mu \subseteq \cl \bigcup_{\ell=1}^\infty \supp\mu_\ell
  .\]
\end{lemma}

The goal of the remaining part of the work is to confirm  Theorem \ref{thm:suppChar} for polyoids. Before we get to the proof, we need to discuss four concepts: projections, cusps, pruning and switching. These will be combined at the end.

\section{Projections}\label{sec:3}
In the following, we assume that $n \ge 1$ and $k \in \N$.
For the proof of Theorem \ref{thm:suppChar} we show two inclusions. For one of these (namely, ``$\subseteq$''), two crucial facts that enable us to carry out the argument are that the touching space (see Definition \ref{def:touching}) of the orthogonal projection of a general convex body $K$ to a linear subspace is the orthogonal projection of the touching space of $K$, which is proved in Lemma \ref{thm:tc3}, and that the orthogonal projection to a subspace of a $k$-polyoid $K$ with generating measure $\mu$   is again a $k$-polyoid for which the projection of $\mu$ is  a generating measure, which is established in Lemma \ref{thm:projVertoid}. Lemmas \ref{thm:tc1}  and \ref{thm:tc2} prepare the proof of Lemma \ref{thm:tc3}. These auxiliary results are treated in the present section. Further ingredients needed to establish the inclusion ``$\subseteq$'' are developed in Section \ref{sec:4}.

\begin{lemma}\label{thm:tc1}
  Let $A \subseteq \R^n$ be a convex set, $W \subseteq \R^n$ a linear subspace and $u \in W \elminus0$.
  Then for all $x \in A$,
  \[
    x \in F(A, u) \iff \pi_W(x) \in F(\pi_W(A), u).
  \]
\end{lemma}
\begin{proof}
  The basic observation is that for all $x \in A$, we have $\left<x, u\right>=\left<\pi_W(x), u\right>$, and hence
    $h_A(u) = h_{\pi_W(A)}(u)$.

  So if $x \in F(A, u)$, then $\left<\pi_W(x), u\right> = \left<x, u\right> = h_A(u) = h_{\pi_W(A)}(u)$ and therefore $\pi_W(x) \in F(\pi_W(A), u)$. Conversely, if $\pi_W(x) \in F(\pi_W(A), u)$, then $\left<x, u\right> = \left<\pi_W(x), u\right> = h_{\pi_W(A)}(u) = h_A(u)$ and hence $x \in F(A, u)$.
\end{proof}

\begin{lemma}\label{thm:tc2}
  Let $K \in \R^n$ be a convex body, $W \subseteq \R^n$ a linear subspace and $u \in W\elminus0$.
  Then $N_W(\pi_{W}(K), F(\pi_W(K), u)) = W \cap N(K, F(K, u))$.
\end{lemma}
\begin{proof}
  By definition of $N_W$, both sides of the equation are subsets of $W$. Moreover, both contain $0$.
  Let $v \in W\elminus0$. Then the claim can be reformulated as
  \[
    F(\pi_W(K), u) \subseteq F(\pi_W(K), v) \iff F(K, u) \subseteq F(K, v)
  .\]

  Let us first assume that $F(\pi_W(K), u) \subseteq F(\pi_W(K), v)$ and let $x \in F(K, u)$. Then by Lemma \ref{thm:tc1}, $\pi_W(x) \in F(\pi_W(K), u)$. By assumption, this implies $\pi_W(x) \in F(\pi_W(K), v)$. Another application of Lemma \ref{thm:tc1} now shows that $x \in F(K, v)$. Therefore, $F(K, u) \subseteq F(K, v)$.

  Now assume $F(K, u) \subseteq F(K, v)$ and let $y \in F(\pi_W(K), u)$. Writing $y = \pi_W(x)$ for some $x \in K$ and applying Lemma \ref{thm:tc1}, we obtain $x \in F(K, u)$. By assumption, this implies $x \in F(K, v)$, and again using Lemma \ref{thm:tc1}, this shows that $y = \pi_W(x) \in F(\pi_W(K), v)$. Therefore, $F(\pi_W(K), u) \subseteq F(\pi_W(K), v)$.
\end{proof}

\begin{lemma}\label{thm:tc3}
  Let $K$ be a convex body, $W \subseteq \R^n$ a linear subspace and $u \in W\elminus0$.
  Then $T_W(\pi_{W}(K), u) = W \cap T(K, u)$ and $\TS_W(\pi_W(K), u) = \pi_W(\TS(K, u))$.
\end{lemma}
\begin{proof}
  By Definition \ref{def:touching}, $T_W(\pi_W(K), u)$ is the unique face of the normal  cone  $N_W(\pi_W(K), F(\pi_W(K), u))$ such that its relative interior contains $u$.
  Similarly, $T(K, u)$ is the unique face of $N(K, F(K, u))$ such that its relative interior contains $u$. We show that $W \cap T(K, u)$ satisfies the definition of $T_W(\pi_W(K), u)$.

  Because $T(K, u)$ is a face of $N(K, F(K, u))$ and by Lemma \ref{thm:tc2},
  \[
    W \cap T(K, u) \text{ is a face of } W \cap N(K, F(K, u)) = N_W(\pi_W(K), F(\pi_W(K), u))
  .\]
  As $(\relint T(K, u)) \cap W$ contains $u$ and $W$ is a linear subspace, $\relint \prn*{W \cap T(K, u)}$ also contains $u$. This proves the first claim.

  For the second claim, observe that $u \in (\relint T(K, u)) \cap W$ implies
  \[
    \vspan (T(K, u) \cap W) = (\vspan T(K, u)) \cap W
  .\]
  Using the first claim, we get
  \[
    \vspan T_W(\pi_W K, u) = \vspan (T(K, u) \cap W) = (\vspan T(K, u)) \cap W
  .\]
  Now we take the orthogonal complement in $W$ and obtain
  \begin{align*}
    \TS_W(\pi_W K, u) &= T_W(\pi_W K, u)^\perp \cap W = (T(K, u)^\perp + W^\perp) \cap W\\
      &= \pi_W T(K, u)^\perp = \pi_W \TS(K, u),
  \end{align*}
  which confirms also the second claim.
\end{proof}

In \cite{HugReichert23+} a $k$-polyoid, for an integer $k\in\N$, was defined as the limit of a sequence of Minkowski sums of $k$-topes, where a $k$-tope is a convex polytope having at most $k$ vertices. Let $\Pt^n_k$ denote the set of $k$-topes in $\R^n$. Furthermore, it was shown in \cite[Thm.~2.8]{HugReichert23+} that a convex body $K\in\K^n$ is a $k$-polyoid if and only if there is a probability measure $\mu$ on $\Pt^n_k$ with compact support such that
\begin{equation}\label{eq:MGM}
h_K(u) = \int h_P(u) \, \mu(\diff P),\quad u \in \R^n.
\end{equation}
Any such (in general non-unique) measure $\mu$ is called a generating measure of the $k$-polyoid $K$.

Let $\varnothing\neq\K_*\subseteq \K^n$ be a Borel set (Borel sets are defined with respect to the topology induced by the Hausdorff metric on $\K^n$). A convex body $K$ in $\R^n$, $n\in\N_0$, for which there  is a probability measure   $\mu$  on $\K_*$ with bounded support such that \eqref{eq:MGM} holds, is called a $\K_*$-\emph{macroid} with generating measure $\mu$. Here the support of $\mu$ is determined with respect to the metric space $\K_*$. It was shown in \cite[Lem.~2.11]{HugReichert23+} that a  $\K_*$-\emph{macroid} with generating measure $\mu$ is the limit of a sequence of Minkowski sums of convex bodies in $\supp\mu$. In the case $\K_*=\Pt^n$, that is, $K$ is a $\Pt^n$-macroid with generating measure $\mu$ on $\Pt^n$, we simply say that $K$ is a macroid with generating measure $\mu$.

\begin{lemma}\label{thm:projVertoid}
  Let $K$ be a macroid (a $k$-polyoid) with generating measure $\mu$,  and let $W \subseteq \R^n$ be a linear subspace.
  Moreover, let $\tilde\pi_W$ be the function that maps the $k$-topes $P \subseteq \R^n$ to the $k$-topes $\pi_W(P) \subseteq W$. Then $\pi_W(K)$ is a macroid (a $k$-polyoid)  with generating measure
  \[
    \mu^W \coloneqq \mu\circ\tilde\pi_W^{-1} \quad \text{and} \quad \tilde\pi_W(\supp\mu)\subseteq\supp\mu^W .
  \]
  If $K \subset \R^n$ is a $k$-polyoid, then  $\tilde\pi_W(\supp\mu)=\supp\mu^W.$
\end{lemma}
\begin{proof} Let $K$ be a macroid (a $k$-polyoid) with generating measure $\mu$.
  For all $u \in W$,
  \[
    h_{\pi_W(K)}(u) = h_K(u) = \int h_P(u) \,\mu(\diff P) = \int h_{\tilde\pi_W(P)}(u) \,\mu(\diff P) = \int h_P(u) \, \mu^W(\diff P)
  .\]
  Moreover,  if $\mu$ is a probability measure  with  bounded (compact) support on polytopes ($k$-topes) in $\R^n$, then $\mu^W$ is a probability measure with  bounded (compact) support on polytopes ($k$-topes) in $W$.

Let $P \in \tilde\pi_W(\supp\mu)$ and $U$ an open neighborhood of $P$ in the space of $k$-topes in $W$. Then there is $Q \in \supp\mu$ such that $\tilde\pi_W(Q) = P$, so that $\tilde\pi_W^{-1}(U)$ is an open neighborhood of $Q$ in $\Pt^n$ (respectively, in $\Pt^n_k$). Therefore,
  \[
    \mu^W(U) = \mu(\tilde\pi_W^{-1}(U)) > 0
  ,\]
  and because this holds for arbitrary $P \in \tilde\pi_W(\supp\mu)$ and neighborhoods $U$ of $P$, it follows that $\tilde\pi_W(\supp\mu) \subseteq \supp\mu^W$.

  Now we assume that $K$ is a $k$-polyoid.
  Because $\tilde\pi_W$ is continuous and $\supp\mu$ is compact, the set $\tilde\pi_W(\supp\mu)$ is compact and hence closed. From
  \[
    \mu^W(\tilde\pi_W(\supp\mu)^\complement) = \mu(\tilde\pi_W^{-1}(\tilde\pi_W(\supp\mu))^\complement) \le \mu((\supp\mu)^\complement) = 0
  \]
  we conclude that $\supp\mu^W \subseteq \tilde\pi_W(\supp\mu)$, and thus    $\supp\mu^W = \tilde\pi_W(\supp\mu)$.
\end{proof}

\begin{remark}
    {\rm Let $C$ be a macroid ($k$-polyoid) with generating measure $\mu$ and $u \in \stdsph$. Recall from \cite[Rem.~2.19]{HugReichert23+} that $F(C, u)$ is a macroid ($k$-polyoid) with generating measure $F_u(\mu) $, which denotes the image measure of $\mu$ under the measurable map $F_u=F(\cdot\,,u)$. In other words,
    \begin{equation}\label{eqsupportset}
        h_{F(C, u)} = \int h_P \,F_u(\mu)(\diff P).
    \end{equation}
  As a consequence, we obtain
    \begin{equation}\label{eqnorcone2}
   \bigcap_{P\in \mathcal{P}(\mu)}N(P,F(P,u))\subseteq  N(C,F(C,u)),
    \end{equation}
    whenever $\mathcal{P}(\mu)\subseteq \Pt^n$ is a measurable set of full $\mu$-measure. For  instance, we can choose $\mathcal{P}(\mu)=\supp \mu$.
  To verify \eqref{eqnorcone2}, let $v\in \bigcap_{P\in \mathcal{P}(\mu)}N(P,F(P,u))$. Then, for each $P\in \mathcal{P}(\mu)$, $F(P,u)\subseteq F(P,v)$, hence $h_{F(P,u)}\le h_{F(P,v)}$. Then \eqref{eqsupportset} yields
  $$
  h_{F(C,u)}=\int h_{F(P,u)}\, \mu(\diff P)
  \le \int h_{F(P,v)}\, \mu(\diff P)= h_{F(C,v)},
  $$
  which shows that $F(C,u)\subseteq F(C,v)$, and therefore $v\in N(C,F(C,u))$.

  A corresponding inclusion for the touching cones does not hold in general, as shown by Example \ref{ex:prune}, which is in contrast to the case of finite Minkowski sums (see \cite[Thm.~2.2.1 (a)]{Schneider}).

  There is a partial converse to \eqref{eqnorcone2}. Let $u\in\mathbb{S}^{n-1}$ be fixed and let $s(L)$ denote the Steiner point of $L\in\K^n$. Recall from \cite[(1.34)]{Schneider} that $s(L)\in \relint L$. Fubini's theorem yields
  $$
  s(F(C,u))=\int s(F(P,u))\, \mu(\diff P),
  $$
  (compare \cite[Rem.~2.14]{HugReichert23+}), and therefore
  $$
  h_{C-s(F(C,u))}=\int h_{P-s(F(P,u))}\, \mu(\diff P).
  $$
  All support functions in this equation are nonnegative. If $v\in N(C,F(C,u))$, then
  $ h_{C-s(F(C,u))}(v)=0$, and hence
 $h_{P-s(F(P,u))}(v)=0$ for $\mu$-almost all $P\in\Pt^n$. This shows that
 $v\in N(P,F(P,u))$ for $\mu$-almost all $P\in\Pt^n$, that is, there is a measurable set $\mathcal{P}_{u,v}(\mu)$ of full $\mu$-measure such that $v\in N(P,F(P,u))$ for all $P\in \mathcal{P}_{u,v}(\mu)$. Let $D_u$ be a countable dense subset of $ N(C,F(C,u))$ and set
 $\mathcal{P}_u(\mu):=\cap_{v\in D_u}\mathcal{P}_{u,v}(\mu)$. Then $\mathcal{P}_{u}(\mu)$ is a measurable set that has full $\mu$-measure and
  $$
   N(C,F(C,u))\subseteq
   \bigcap_{P\in \mathcal{P}_u(\mu)}N(P,F(P,u)).
 $$
Together with \eqref{eqnorcone2} we obtain
  $$
   N(C,F(C,u))=
   \bigcap_{P\in \mathcal{P}_u(\mu)}N(P,F(P,u)).
 $$
  }
\end{remark}

\section{Cusps}\label{sec:4}

The proof of Theorem \ref{thm:suppChar} relies on the assumption that the convex bodies in question are  polyoids. In fact, one inclusion holds for the larger class of macroids. For this reason, the results in this section are provided for the class of macroids or  for general convex bodies.    The following results about \emph{cusps} describe what it means that the touching space of a convex body $K$ (a macroid $K$ with generating measure $\mu$) is $0$-dimensional, in terms of the polytopes in the support of $\mu$. One might hope that $\TS(K, u) = \set*{0}$ if and only if the same holds for all $P \in \supp\mu$, but this turns out to be false (the \enquote{only if} statement is true though, as follows from Lemmas \ref{thm:tc4} and \ref{thm:tc5}). Cusps can be thought of as an attempt to quantify how far a convex body is from having a non-trivial touching space. Intuitively, Lemmas \ref{thm:tc4} and \ref{thm:tc5} show that $\TS(K, u)$ is trivial if and only if the $k$-topes in $\supp\mu$ keep a minimum distance from having a non-trivial touching space.
Lemmas \ref{thm:tc4} and \ref{thm:tc5}  will be employed in the crucial Witness Lemma \ref{thm:pruningLemma}.

\begin{definition}
{\rm
  For all $u \in \stdsph$ and $c > 0$, define a cone with apex at $0$,
  \[
    \ccone_c(u) \coloneqq \set{ x \in \R^n \given \left<x, u\right> \le -c \norm{x} }
  .\]
  Let $K \subseteq  \R^n$ be a convex body, $u \in \stdsph$ and $c > 0$. Then $K$ is said to \emph{have a $c$-cusp in direction $u \in \stdsph$} if there is some $x \in K$ such that $K \subseteq x + \ccone_c(u)$.}
\end{definition}

Note that $ \ccone_c(u)=\{0\}$ if $c>1$ and $\ccone_1(u)=-[0,\infty) u$; the cone $ \ccone_c(u)$ is getting smaller as $c\in (0,1]$ is getting larger. In particular, if $K$ has a $c$-cusp in direction $u$, then it also has a $c^\prime$-cusp in direction $u$ for $0<c^\prime<c$.

\begin{lemma}
\label{thm:locLin}
  Let $K \in \K^n$ be a convex body, $u \in \stdsph$ and $c > 0$.
  Then the following are equivalent:
  \begin{enumerate}[{\rm (a)}]
    \item $K$ has a $c$-cusp in direction $u$.
    \item $h_K$ is linear on $U(u, c) \coloneqq c\stdb + u$.
  \end{enumerate}
\end{lemma}
\begin{proof}
  The statement is invariant under translations.

\noindent
  \enquote{(a) $\implies$ (b)}: Assume that there is some $x \in K$ with $K \subseteq x + \ccone_c(u)$. Translating $K$, we can arrange that $x = 0$.
  Then the Cauchy--Schwarz inequality shows that for all $v \in U(u, c)$ and $y \in K\subseteq \ccone_c(u)$,
  \[
    \left<y, v\right> \le \left<y, u\right> + \norm{u - v}\norm{y} \le (-c + \norm{u - v}) \norm{y} \le 0 = \left<x, v\right>
  .\]
  So $h_K(v) = 0$ for $v\in U(u, c)$.

  \medskip

\noindent
  \enquote{(b) $\implies$ (a)}: Assume that there is some $x \in \R^n$ such that $h_K = \left<x, \cdot\right>$ on $U(u, c)$. Translating $K$ by $-x$, we can arrange that $x = 0$. Then for all $y \in K \elminus0$,
  \[
    \left<y, u\right> = \left\langle y, u + \frac{c}{\norm{y}} y\right\rangle - c{\norm{y}}   \le h_K\left( u + \frac{c}{\norm{y}} y \right) -  {c}\norm{y}  = -c\norm{y}
  .\]
  So $K \subseteq \ccone_c(u)$ (remembering $0 \in \ccone_c(u)$). Moreover, $h_K'(u; \cdot) = 0$ because $U(u,c)$ is a neighborhood of $u$ where $h_K \equiv 0$. With \cite[Thm.~1.7.2]{Schneider} it follows that
  \[
    h_{F(K, u)} = h_K'(u; \cdot) = 0 = h_{\set*{0}}
  ,\]
  proving that $0 \in \set*{0} = F(K, u) \subseteq K$. So $K \subseteq \ccone_c(u)$ and $0 \in K$.
\end{proof}

Next we use Lemma \ref{thm:locLin}  to characterize the situation when the touching space is trivial.

\begin{lemma}
\label{thm:tc4}
  Let $K \in \K^n$ be a convex body, and let $u \in \mathbb{S}^{n-1}$. Then the following are equivalent.
  \begin{enumerate}[{\rm (a)}]
    \item $\TS(K, u) = \set*{0}$.
    \item There is some $c > 0$ such that $K$ has a $c$-cusp in direction $u$.
  \end{enumerate}
\end{lemma}

\begin{proof}
\enquote{(a) $\implies$ (b)}:
  Assume that $\TS(K, u) = \set*{0}$. Then $u \in \absint N(K, F(K, u))$. So there is $c > 0$ such that $U(u, c) = \set*{u} + c\stdb \subseteq N(K, F(K, u))$. Choosing $x \in F(K, u)$, it follows that $h_K = \left<x, \cdot\right>$ on $U(u, c) \subseteq N(K, F(K, u))$. Then by Lemma \ref{thm:locLin}, $K$ has a $c$-cusp in direction $u$.

  \medskip

\noindent
\enquote{(b) $\implies$ (a)}:
  Assume that $K$ has a $c$-cusp in direction $u$ for some $c > 0$. Then by Lemma \ref{thm:locLin}, there is $x \in \R^n$ such that $h_K = \left<x, \cdot\right>$ on $U(u, c) = \set*{u} + c\stdb$. By \cite[Thm.~1.7.2]{Schneider}, all $v \in \absint U(u, c)$ satisfy $h_{F(K, v)} = h_K'(v; \cdot) = \left<x, \cdot\right>$, so that $F(K, v) = \set*{x} = F(K, u)$. Hence, $\absint(\set*{u} + c\stdb) \subseteq N(K, F(K, u))$, showing that $u \in \absint N(K, F(K, u))$ and $\TS(K, u) = \set*{0}$.
\end{proof}

In the following we need to understand how the local linearity of the support function of a macroid is related to the local linearity of the support functions of the polytopes in the support of a generating measure of the macroid. This relation is given in Lemma \ref{thm:locLin2}, which we prepare by two simple lemmas.

The first lemma is well-known, but we state it for easier reference. The proof of the second lemma is included, since it is crucial for the proof of Lemma \ref{thm:locLin2}.

\begin{lemma}\label{thm:suppAff}
  Let $A \subseteq \R^n$ be a convex set, $f \colon A \to \R$ a convex function and $a \in \relint A$. Then there is $u \in \pspan A$ such that
  \[
    f(x) \ge \left<x - a, u\right> + f(a) \quad \text{for all $x \in A$}
  .\]
\end{lemma}

\begin{lemma}
\label{thm:locLin3}
  Let $A \subseteq \R^n$ be a convex set, and let $f\colon \R^n \to \R$ be positively $1$-homogeneous. Then the following are equivalent.
  \begin{enumerate}[{\rm (a)}]
    \item\label{it:conc1} $f$ is linear on $A$ (i.e.\ agrees on $A$ with a function $x \mapsto \left<x, u\right>$, where $u \in \R^n$).
    \item\label{it:conc2} $f$ is affine on $A$ (i.e.\ agrees on $A$ with a function $x \mapsto \left<x, u\right> + c$, where $u \in \R^n$ and $c \in \R$).
    \item\label{it:conc3} $f$ is convex and concave on $A$.
  \end{enumerate}
\end{lemma}
\begin{proof}
  (\ref{it:conc1}) implies (\ref{it:conc2}) and (\ref{it:conc2}) implies (\ref{it:conc3}). Without loss of generality, $A$ is nonempty.

  \medskip

\noindent
  \enquote{(\ref{it:conc2}) $\implies$ (\ref{it:conc1})}: Assume that there are $u \in \R^n$ and $c \in \R$ such that
  \[
    f(x) = \left<x, u\right> + c \quad \text{for $x \in A$}
  .\]
  Let $E$ be the affine span of $A$. If $0 \in E$, then choose $x \in \relint A$. There is $\lambda \in (0, 1)$ such that $\lambda x \in A$, so that we obtain
  \[
    \lambda\left<x, u\right> + c = f(\lambda x) = \lambda f(x) = \lambda\left<x, u\right> + \lambda c \implies c = \lambda c \implies c = 0
  .\]
  If $0 \notin E$, then $E \cap \pspan A = \varnothing$. Choose $a \in A$. Then $a \notin \pspan A = \prn*{\pspan A}^{\perp\perp}$ and there is $v \in \prn*{\pspan A}^\perp$ such that $\left<a, v\right> \ne 0$. Also observe that if $x \in A$, then $x - a \in \pspan A$, so that $\left<x, v\right> = \left<a, v\right>$. So
  \[
    f(x) = \left<x, u\right> + c = \left<x, u\right> + c\frac{\left<x, v\right>}{\left<a, v\right>} = \left<x, u + c\frac{v}{\left<a, v\right>}\right> \quad \text{for all $x \in A$}
  .\]

\noindent
  \enquote{(\ref{it:conc3}) $\implies$ (\ref{it:conc2})}:
  Let $a \in \relint A$. By convexity of $f$ and Lemma \ref{thm:suppAff}, there is $u \in \pspan A$ such that
  \[
    f(x) \ge \left<x - a, u\right> + f(a) \quad \text{for all $x \in A$}
  .\]
  By concavity of $f$ and Lemma \ref{thm:suppAff} applied to $-f$, there is $v \in \pspan A$ such that
  \[
    f(x) \le \left<x - a, v\right> + f(a) \quad \text{for all $x \in A$}
  .\]
  Hence, $\left<x - a, v - u\right> \ge 0$ for all $x \in A$. Because $a \in \relint A$ and $u, v \in \pspan A$, this shows that $v = u$ and so
  \[
    f(x) = \left<x - a, u\right> + f(a) = \left<x, u\right> - \left<a, u\right> + f(a) \quad \text{for all $x \in A$}
  ,\]
  which completes the proof.
\end{proof}

\begin{lemma}
\label{thm:locLin2}
Let $\varnothing\neq\K_*\subseteq \K^n$ be a Borel set. Let $K\in\K^n$ be a $\K_*$-macroid with generating measure  $\mu$ on $\K_*$, and let $A \subseteq \R^n$ be convex. Then $h_K$ is linear on $A$ if and only if $h_P$ is linear on $A$ for all $P \in \supp\mu$.
\end{lemma}
\begin{proof}
  Every support function of a convex body is convex and positively $1$-homogeneous. So by Lemma \ref{thm:locLin3}, it is linear on $A$ if and only if it is concave on $A$.

\noindent
  \enquote{$\implies$}: Assume that there is $P \in \supp\mu$ such that $h_P$ is not concave on $A$. Then there are an open neighborhood $U$ of $P$, $\lambda \in (0, 1)$ and $y, z \in A$ such that for all $Q \in U$,
  \[
    h_Q(\lambda y + (1 - \lambda) z) < \lambda h_Q(y) + (1 - \lambda) h_Q(z)
  .\]
  On the other hand, for all $Q \in U^\complement$, convexity of $h_Q$ implies
  \[
    h_Q(\lambda y + (1 - \lambda) z) \le \lambda h_Q(y) + (1 - \lambda) h_Q(z)
  .\]
  Since $\mu(U) > 0$, we thus obtain from \eqref{eq:MGM} that
  \[
    h_K(\lambda y + (1 - \lambda) z) < \lambda h_K(y) + (1 - \lambda) h_K(z)
  .\]
  Therefore, $h_K$ is not concave on $A$.

\medskip

\noindent
  \enquote{$\impliedby$}: Assume that $h_K$ is not concave on $A$. Then there are $\lambda \in (0, 1)$ and $y, z \in A$ such that
  \[
    h_K(\lambda y + (1 - \lambda) z) < \lambda h_K(y) + (1 - \lambda) h_K(z)
  .\]
  In particular, there is at least one $P \in \supp\mu$ with
  \[
    h_P(\lambda y + (1 - \lambda) z) < \lambda h_P(y) + (1 - \lambda) h_P(z)
  .\]
  Therefore, $h_P$ is not concave on $A$.
\end{proof}

\begin{lemma}
\label{thm:tc5}
 Let $\varnothing\neq\K_*\subseteq \K^n$ be a Borel set. Let $K\in\K^n$ be a $\K_*$-macroid with generating measure  $\mu$ on $\K_*$. Let $u \in \mathbb{S}^{n-1}$ and $c > 0$. Then the following are equivalent.
  \begin{enumerate}[{\rm (a)}]
    \item $K$ has a $c$-cusp in direction $u$.
    \item Every $P \in \supp\mu$ has a $c$-cusp in direction $u$.
  \end{enumerate}
\end{lemma}
\begin{proof}
  By Lemma \ref{thm:locLin}, $K$ has a $c$-cusp in direction $u$ if and only if $h_K$ is linear on $U(u, c)$. By Lemma \ref{thm:locLin2}, this is equivalent to $h_P$ being linear on $U(u, c)$ for all $P \in \supp\mu$. Again by Lemma \ref{thm:locLin}, this in turn is equivalent to $P$ having a $c$-cusp in direction $u$ for all $P \in \supp\mu$.
\end{proof}

As a consequence of Lemmas \ref{thm:tc4} and \ref{thm:tc5}, we obtain the following corollary.

\begin{corollary}\label{cor:cusp}
 Let $\varnothing\neq\K_*\subseteq \K^n$ be a Borel set.  Let $K\in\K^n$ be a $\K_*$-macroid  with generating measure  $\mu$ on $\K_*$ and let   $u \in \mathbb{S}^{n-1}$. Then the following statements are equivalent:
\begin{enumerate}[{\rm (a)}]
    \item $\TS(K,u)\neq \{0\}$.
\item For each $c > 0$ there exists some $P \in\supp \mu$ that does not have a $c$-cusp in direction $u$.
\end{enumerate}
\end{corollary}

We denote by $\K^n_{sm}$ the set of all smooth convex bodies. Since the complement of
$\K^n_{sm}$ is a countable union of closed sets, $\K^n_{sm}$ is measurable. It follows from \cite[Thm.~2.2.1 (a)]{Schneider} that a finite Minkowski sum of convex bodies, one of which is smooth, is smooth again. In other words, if the sum is not smooth, then none of the summands is smooth.  Next we show that this fact extends to macroids. In particular, there is no point in considering $\K^n_{sm}$-macroids.

\begin{corollary}\label{corsmooth}
Suppose that $K$ is a $\K_*$-macroid with generating measure $\mu$ and $K$ is not smooth. Then none of the $L\in\supp \mu$ is smooth.
\end{corollary}

\begin{proof}
    If $K$ is not smooth, then there is a convex cone $A$ with $\dim A\ge 2$ such that $h_K$ is linear on $A$. By Lemma \ref{thm:tc5}, $h_P$ is linear on $A$, for each $P\in \supp\mu$. But then $P$ is not smooth, for each $P\in\supp\mu$.
\end{proof}

\section{Pruning}\label{sec:5}

This section develops a technique that is only relevant for proving one of the two inclusions on which the characterization Theorem \ref{thm:suppChar}   is based: Let $\collc = (C_1, \ldots, C_{n-1})$ be a tuple of $k$-polyoids with generating measures $\mu_1, \ldots, \mu_{n-1}$. If $u \in \ext \collc$, then we have to show that $u \in \supp\Su(\collc)$.

Because this is the most difficult aspect of \nref{thm:suppChar}, we begin with some examples. The first example introduces the idea of a \enquote{witness polytope} that is used to prove that some normal vector is in the support of a mixed area measure. The other two examples exemplify how to find \enquote{witness polytopes} in more complicated situations using \emph{pruning}, the method developed in this section.

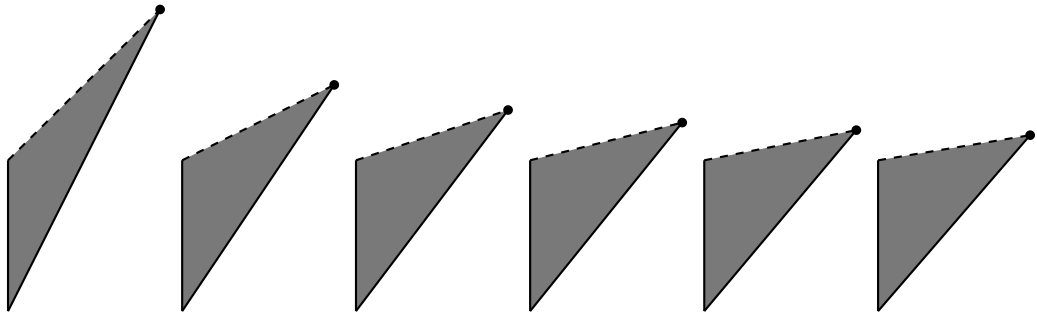
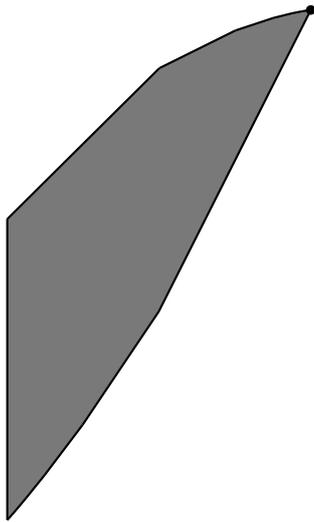
\begin{figure}
    \centering
    \begin{subfigure}[b]{\textwidth}
      \centering
      \begin{tikzpicture}%
	[scale=2.000000,
	back/.style={loosely dotted, thin},
	edge/.style={color=black, thick},
	facet/.style={fill=gray!95!black,fill opacity=0.800000},
	vertex/.style={inner sep=1pt,circle,draw=black!25!black,fill=black!75!black,thick,anchor=base}]
%
%
\coordinate (0.00000, 0.00000) at (0.00000, 0.00000);
\coordinate (0.00000, 1.00000) at (0.00000, 1.00000);
\coordinate (1.00000, 2.00000) at (1.00000, 2.00000);
\fill[facet] (1.00000, 2.00000) -- (0.00000, 0.00000) -- (0.00000, 1.00000) -- cycle {};
\draw[edge] (0.00000, 0.00000) -- (0.00000, 1.00000);
\draw[edge] (0.00000, 0.00000) -- (1.00000, 2.00000);
\draw[edge,dashed] (0.00000, 1.00000) -- (1.00000, 2.00000);
\node[vertex] at (1.00000, 2.00000)     {};
\end{tikzpicture}
      \begin{tikzpicture}%
	[scale=2.000000,
	back/.style={loosely dotted, thin},
	edge/.style={color=black, thick},
	facet/.style={fill=gray!95!black,fill opacity=0.800000},
	vertex/.style={inner sep=1pt,circle,draw=black!25!black,fill=black!75!black,thick,anchor=base}]
%
%
\coordinate (0.00000, 0.00000) at (0.00000, 0.00000);
\coordinate (0.00000, 1.00000) at (0.00000, 1.00000);
\coordinate (1.00000, 1.50000) at (1.00000, 1.50000);
\fill[facet] (1.00000, 1.50000) -- (0.00000, 0.00000) -- (0.00000, 1.00000) -- cycle {};
\draw[edge] (0.00000, 0.00000) -- (0.00000, 1.00000);
\draw[edge] (0.00000, 0.00000) -- (1.00000, 1.50000);
\draw[edge,dashed] (0.00000, 1.00000) -- (1.00000, 1.50000);
\node[vertex] at (1.00000, 1.50000)     {};
\end{tikzpicture}
      \begin{tikzpicture}%
	[scale=2.000000,
	back/.style={loosely dotted, thin},
	edge/.style={color=black, thick},
	facet/.style={fill=gray!95!black,fill opacity=0.800000},
	vertex/.style={inner sep=1pt,circle,draw=black!25!black,fill=black!75!black,thick,anchor=base}]
%
%
\coordinate (0.00000, 0.00000) at (0.00000, 0.00000);
\coordinate (0.00000, 1.00000) at (0.00000, 1.00000);
\coordinate (1.00000, 1.33333) at (1.00000, 1.33333);
\fill[facet] (1.00000, 1.33333) -- (0.00000, 0.00000) -- (0.00000, 1.00000) -- cycle {};
\draw[edge] (0.00000, 0.00000) -- (0.00000, 1.00000);
\draw[edge] (0.00000, 0.00000) -- (1.00000, 1.33333);
\draw[edge,dashed] (0.00000, 1.00000) -- (1.00000, 1.33333);
\node[vertex] at (1.00000, 1.33333)     {};
\end{tikzpicture}
      \begin{tikzpicture}%
	[scale=2.000000,
	back/.style={loosely dotted, thin},
	edge/.style={color=black, thick},
	facet/.style={fill=gray!95!black,fill opacity=0.800000},
	vertex/.style={inner sep=1pt,circle,draw=black!25!black,fill=black!75!black,thick,anchor=base}]
%
%
\coordinate (0.00000, 0.00000) at (0.00000, 0.00000);
\coordinate (0.00000, 1.00000) at (0.00000, 1.00000);
\coordinate (1.00000, 1.25000) at (1.00000, 1.25000);
\fill[facet] (1.00000, 1.25000) -- (0.00000, 0.00000) -- (0.00000, 1.00000) -- cycle {};
\draw[edge] (0.00000, 0.00000) -- (0.00000, 1.00000);
\draw[edge] (0.00000, 0.00000) -- (1.00000, 1.25000);
\draw[edge,dashed] (0.00000, 1.00000) -- (1.00000, 1.25000);
\node[vertex] at (1.00000, 1.25000)     {};
\end{tikzpicture}
      \begin{tikzpicture}%
	[scale=2.000000,
	back/.style={loosely dotted, thin},
	edge/.style={color=black, thick},
	facet/.style={fill=gray!95!black,fill opacity=0.800000},
	vertex/.style={inner sep=1pt,circle,draw=black!25!black,fill=black!75!black,thick,anchor=base}]
%
%
\coordinate (0.00000, 0.00000) at (0.00000, 0.00000);
\coordinate (0.00000, 1.00000) at (0.00000, 1.00000);
\coordinate (1.00000, 1.20000) at (1.00000, 1.20000);
\fill[facet] (1.00000, 1.20000) -- (0.00000, 0.00000) -- (0.00000, 1.00000) -- cycle {};
\draw[edge] (0.00000, 0.00000) -- (0.00000, 1.00000);
\draw[edge] (0.00000, 0.00000) -- (1.00000, 1.20000);
\draw[edge,dashed] (0.00000, 1.00000) -- (1.00000, 1.20000);
\node[vertex] at (1.00000, 1.20000)     {};
\end{tikzpicture}
      \begin{tikzpicture}%
	[scale=2.000000,
	back/.style={loosely dotted, thin},
	edge/.style={color=black, thick},
	facet/.style={fill=gray!95!black,fill opacity=0.800000},
	vertex/.style={inner sep=1pt,circle,draw=black!25!black,fill=black!75!black,thick,anchor=base}]
%
%
\coordinate (0.00000, 0.00000) at (0.00000, 0.00000);
\coordinate (0.00000, 1.00000) at (0.00000, 1.00000);
\coordinate (1.00000, 1.16667) at (1.00000, 1.16667);
\fill[facet] (1.00000, 1.16667) -- (0.00000, 0.00000) -- (0.00000, 1.00000) -- cycle {};
\draw[edge] (0.00000, 0.00000) -- (0.00000, 1.00000);
\draw[edge] (0.00000, 0.00000) -- (1.00000, 1.16667);
\draw[edge,dashed] (0.00000, 1.00000) -- (1.00000, 1.16667);
\node[vertex] at (1.00000, 1.16667)     {};
\end{tikzpicture}
      \caption{The first six triangles $C^\tidx\ell$ from \cref{ex:simpleWitness}. The faces spanned by the second and third vertex are dashed; they converge to a nondegenerate $e_2$-face of the limit triangle. The support sets in direction $e_2$ consist of a single point marked black.}
    \end{subfigure}
    \begin{subfigure}[b]{0.95\textwidth}
      \centering
    \begin{tikzpicture}%
	[scale=2.000000,
	back/.style={loosely dotted, thin},
	edge/.style={color=black, thick},
	facet/.style={fill=gray!95!black,fill opacity=0.800000},
	vertex/.style={inner sep=1pt,circle,draw=black!25!black,fill=black!75!black,thick,anchor=base}]
%
%
\coordinate (0.00000, 0.00000) at (0.00000, 0.00000);
\coordinate (0.00391, 0.00434) at (0.00391, 0.00434);
\coordinate (0.01172, 0.01313) at (0.01172, 0.01313);
\coordinate (0.02734, 0.03099) at (0.02734, 0.03099);
\coordinate (0.05859, 0.06744) at (0.05859, 0.06744);
\coordinate (0.12109, 0.14244) at (0.12109, 0.14244);
\coordinate (0.24609, 0.29869) at (0.24609, 0.29869);
\coordinate (0.49609, 0.63203) at (0.49609, 0.63203);
\coordinate (0.99609, 1.38203) at (0.99609, 1.38203);
\coordinate (0.00000, 1.99609) at (0.00000, 1.99609);
\coordinate (1.00000, 2.99609) at (1.00000, 2.99609);
\coordinate (1.50000, 3.24609) at (1.50000, 3.24609);
\coordinate (1.75000, 3.32943) at (1.75000, 3.32943);
\coordinate (1.87500, 3.36068) at (1.87500, 3.36068);
\coordinate (1.93750, 3.37318) at (1.93750, 3.37318);
\coordinate (1.96875, 3.37839) at (1.96875, 3.37839);
\coordinate (1.98438, 3.38062) at (1.98438, 3.38062);
\coordinate (1.99219, 3.38159) at (1.99219, 3.38159);
\coordinate (1.99609, 3.38203) at (1.99609, 3.38203);
\fill[facet] (1.99609, 3.38203) -- (0.99609, 1.38203) -- (0.49609, 0.63203) -- (0.24609, 0.29869) -- (0.12109, 0.14244) -- (0.05859, 0.06744) -- (0.02734, 0.03099) -- (0.01172, 0.01313) -- (0.00391, 0.00434) -- (0.00000, 0.00000) -- (0.00000, 1.99609) -- (1.00000, 2.99609) -- (1.50000, 3.24609) -- (1.75000, 3.32943) -- (1.87500, 3.36068) -- (1.93750, 3.37318) -- (1.96875, 3.37839) -- (1.98438, 3.38062) -- (1.99219, 3.38159) -- cycle {};
\draw[edge] (0.00000, 0.00000) -- (0.00391, 0.00434);
\draw[edge] (0.00000, 0.00000) -- (0.00000, 1.99609);
\draw[edge] (0.00391, 0.00434) -- (0.01172, 0.01313);
\draw[edge] (0.01172, 0.01313) -- (0.02734, 0.03099);
\draw[edge] (0.02734, 0.03099) -- (0.05859, 0.06744);
\draw[edge] (0.05859, 0.06744) -- (0.12109, 0.14244);
\draw[edge] (0.12109, 0.14244) -- (0.24609, 0.29869);
\draw[edge] (0.24609, 0.29869) -- (0.49609, 0.63203);
\draw[edge] (0.49609, 0.63203) -- (0.99609, 1.38203);
\draw[edge] (0.99609, 1.38203) -- (1.99609, 3.38203);
\draw[edge] (0.00000, 1.99609) -- (1.00000, 2.99609);
\draw[edge] (1.00000, 2.99609) -- (1.50000, 3.24609);
\draw[edge] (1.50000, 3.24609) -- (1.75000, 3.32943);
\draw[edge] (1.75000, 3.32943) -- (1.87500, 3.36068);
\draw[edge] (1.87500, 3.36068) -- (1.93750, 3.37318);
\draw[edge] (1.93750, 3.37318) -- (1.96875, 3.37839);
\draw[edge] (1.96875, 3.37839) -- (1.98438, 3.38062);
\draw[edge] (1.98438, 3.38062) -- (1.99219, 3.38159);
\draw[edge] (1.99219, 3.38159) -- (1.99609, 3.38203);
\node[vertex] at (1.99609, 3.38203)     {};
\end{tikzpicture}
    \caption{An approximation to the triangle body $C$ obtained by weighted summation of the triangles $C^\tidx\ell$. The support set in direction $e_2$ consists of a single point marked black.}
    \end{subfigure}
    \caption{The situation of Example \ref{ex:simpleWitness}}
    \label{fig:ex1}
  \end{figure}

\begin{example}[A witness polytope]\label{ex:simpleWitness} {\rm
  Let $n=2$. Let $(e_1,e_2)$ be the standard orthonormal basis of $\R^2$.
  Let
  \[
    C^\tidx\ell \coloneqq \conv\set*{0, e_2, e_1+ (1+ {\ell}^{-1})e_2},\quad \ell\in\N,
  \]
  and define the triangle body (i.e., the $3$-polyoid)
  \[
    C \coloneqq \sum_{\ell=1}^\infty 2^{-\ell} C^\tidx\ell
  \]
  with generating measure
  \[
    \mu \coloneqq \sum_{\ell=1}^\infty  2^{-\ell} \delta_{C^\tidx\ell}
  ;\]
  see Figure \ref{fig:ex1} for an illustration.
  The sequence $(C^\tidx\ell)_\ell$ converges to the triangle
  \[
    K \coloneqq \conv\set*{0, e_2, e_1+e_2}
  \]
  and so
  \[
    \supp\mu = \set*{K, C^\tidx 1, C^\tidx 2, \ldots}
  .\]
  By Corollary \ref{cor:cusp} we find that $\TS(C, e_2) \ne \set*{0}$ because $K \in\supp\mu$ does not have a $c$-cusp in direction $e_2$ for any $c > 0$. Hence, $e_2$ is a $(C)$-extreme normal vector.
  So Theorem \ref{thm:suppChar} predicts $e_2 \in \supp\Su(C)$.

  Indeed, Theorem \ref{thm:suppInt} and $K \in \supp\mu$ show that
  \[
    e_2 \in \supp\Su(K) \subseteq \supp\Su(C)
  .\]
  Alternatively, we could argue that $C^\tidx\ell \to K$ and so Lemma \ref{thm:suppEqOfTendsto} and Theorem \ref{thm:suppInt} yield
  \[
    e_2 \in \supp\Su(K) \subseteq \bigcup_{\ell=1}^\infty \supp\Su(C^\tidx\ell) \subseteq \bigcup_{P\in\supp\mu} \supp\Su(P) = \supp\Su(C)
  .\]
  We have used $K \in \supp\mu$ as a \enquote{witness polytope} to establish $e_2 \in \supp\Su(C)$.}
\end{example}

\begin{figure}
  \centering
  \begin{subfigure}[b]{0.5\textwidth}
    \centering
    \begin{tikzpicture}%
	[scale=2.000000,
	back/.style={loosely dotted, thin},
	edge/.style={color=black, thick},
	facet/.style={fill=gray!95!black,fill opacity=0.800000},
	vertex/.style={inner sep=1pt,circle,draw=black!25!black,fill=black!75!black,thick,anchor=base}]
%
%
\coordinate (0.00000, -1.00000) at (0.00000, -1.00000);
\coordinate (0.00000, 0.00000) at (0.00000, 0.00000);
\coordinate (-1.00000, -1.00000) at (-1.00000, -1.00000);
\fill[facet] (-1.00000, -1.00000) -- (0.00000, -1.00000) -- (0.00000, 0.00000) -- cycle {};
\draw[edge] (0.00000, -1.00000) -- (0.00000, 0.00000);
\draw[edge] (0.00000, -1.00000) -- (-1.00000, -1.00000);
\draw[edge,dashed] (0.00000, 0.00000) -- (-1.00000, -1.00000);
\node[vertex] at (0.00000, 0.00000)     {};
\end{tikzpicture}
    \begin{tikzpicture}%
	[scale=2.000000,
	back/.style={loosely dotted, thin},
	edge/.style={color=black, thick},
	facet/.style={fill=gray!95!black,fill opacity=0.800000},
	vertex/.style={inner sep=1pt,circle,draw=black!25!black,fill=black!75!black,thick,anchor=base}]
%
%
\coordinate (0.00000, -1.00000) at (0.00000, -1.00000);
\coordinate (0.00000, 0.00000) at (0.00000, 0.00000);
\coordinate (-0.50000, -0.25000) at (-0.50000, -0.25000);
\fill[facet] (-0.50000, -0.25000) -- (0.00000, -1.00000) -- (0.00000, 0.00000) -- cycle {};
\draw[edge] (0.00000, -1.00000) -- (0.00000, 0.00000);
\draw[edge] (0.00000, -1.00000) -- (-0.50000, -0.25000);
\draw[edge,dashed] (0.00000, 0.00000) -- (-0.50000, -0.25000);
\node[vertex] at (0.00000, 0.00000)     {};
\end{tikzpicture}
    \begin{tikzpicture}%
	[scale=2.000000,
	back/.style={loosely dotted, thin},
	edge/.style={color=black, thick},
	facet/.style={fill=gray!95!black,fill opacity=0.800000},
	vertex/.style={inner sep=1pt,circle,draw=black!25!black,fill=black!75!black,thick,anchor=base}]
%
%
\coordinate (0.00000, -1.00000) at (0.00000, -1.00000);
\coordinate (0.00000, 0.00000) at (0.00000, 0.00000);
\coordinate (-0.33333, -0.11111) at (-0.33333, -0.11111);
\fill[facet] (-0.33333, -0.11111) -- (0.00000, -1.00000) -- (0.00000, 0.00000) -- cycle {};
\draw[edge] (0.00000, -1.00000) -- (0.00000, 0.00000);
\draw[edge] (0.00000, -1.00000) -- (-0.33333, -0.11111);
\draw[edge,dashed] (0.00000, 0.00000) -- (-0.33333, -0.11111);
\node[vertex] at (0.00000, 0.00000)     {};
\end{tikzpicture}
    \begin{tikzpicture}%
	[scale=2.000000,
	back/.style={loosely dotted, thin},
	edge/.style={color=black, thick},
	facet/.style={fill=gray!95!black,fill opacity=0.800000},
	vertex/.style={inner sep=1pt,circle,draw=black!25!black,fill=black!75!black,thick,anchor=base}]
%
%
\coordinate (0.00000, -1.00000) at (0.00000, -1.00000);
\coordinate (0.00000, 0.00000) at (0.00000, 0.00000);
\coordinate (-0.25000, -0.06250) at (-0.25000, -0.06250);
\fill[facet] (-0.25000, -0.06250) -- (0.00000, -1.00000) -- (0.00000, 0.00000) -- cycle {};
\draw[edge] (0.00000, -1.00000) -- (0.00000, 0.00000);
\draw[edge] (0.00000, -1.00000) -- (-0.25000, -0.06250);
\draw[edge,dashed] (0.00000, 0.00000) -- (-0.25000, -0.06250);
\node[vertex] at (0.00000, 0.00000)     {};
\end{tikzpicture}
    \begin{tikzpicture}%
	[scale=2.000000,
	back/.style={loosely dotted, thin},
	edge/.style={color=black, thick},
	facet/.style={fill=gray!95!black,fill opacity=0.800000},
	vertex/.style={inner sep=1pt,circle,draw=black!25!black,fill=black!75!black,thick,anchor=base}]
%
%
\coordinate (0.00000, -1.00000) at (0.00000, -1.00000);
\coordinate (0.00000, 0.00000) at (0.00000, 0.00000);
\coordinate (-0.20000, -0.04000) at (-0.20000, -0.04000);
\fill[facet] (-0.20000, -0.04000) -- (0.00000, -1.00000) -- (0.00000, 0.00000) -- cycle {};
\draw[edge] (0.00000, -1.00000) -- (0.00000, 0.00000);
\draw[edge] (0.00000, -1.00000) -- (-0.20000, -0.04000);
\draw[edge,dashed] (0.00000, 0.00000) -- (-0.20000, -0.04000);
\node[vertex] at (0.00000, 0.00000)     {};
\end{tikzpicture}
    \begin{tikzpicture}%
	[scale=2.000000,
	back/.style={loosely dotted, thin},
	edge/.style={color=black, thick},
	facet/.style={fill=gray!95!black,fill opacity=0.800000},
	vertex/.style={inner sep=1pt,circle,draw=black!25!black,fill=black!75!black,thick,anchor=base}]
%
%
\coordinate (0.00000, -1.00000) at (0.00000, -1.00000);
\coordinate (0.00000, 0.00000) at (0.00000, 0.00000);
\coordinate (-0.16667, -0.02778) at (-0.16667, -0.02778);
\fill[facet] (-0.16667, -0.02778) -- (0.00000, -1.00000) -- (0.00000, 0.00000) -- cycle {};
\draw[edge] (0.00000, -1.00000) -- (0.00000, 0.00000);
\draw[edge] (0.00000, -1.00000) -- (-0.16667, -0.02778);
\draw[edge,dashed] (0.00000, 0.00000) -- (-0.16667, -0.02778);
\node[vertex] at (0.00000, 0.00000)     {};
\end{tikzpicture}
    \caption{The first six triangles $C^\tidx\ell$ from Example \ref{ex:prune}. The faces spanned by the second and third vertex are dashed. The support sets in direction $e_2$ consist of a single point marked black.}
  \end{subfigure}\hspace*{5mm}
  \begin{subfigure}[b]{0.4\textwidth}
    \centering
  \begin{tikzpicture}%
	[scale=2.000000,
	back/.style={loosely dotted, thin},
	edge/.style={color=black, thick},
	facet/.style={fill=gray!95!black,fill opacity=0.800000},
	vertex/.style={inner sep=1pt,circle,draw=black!25!black,fill=black!75!black,thick,anchor=base}]
%
%
\coordinate (0.00000, -1.99994) at (0.00000, -1.99994);
\coordinate (0.00000, 0.00000) at (0.00000, 0.00000);
\coordinate (-0.00072, 0.00000) at (-0.00072, 0.00000);
\coordinate (-0.00169, -0.00012) at (-0.00169, -0.00012);
\coordinate (-0.00393, -0.00044) at (-0.00393, -0.00044);
\coordinate (-0.00914, -0.00131) at (-0.00914, -0.00131);
\coordinate (-0.02164, -0.00381) at (-0.02164, -0.00381);
\coordinate (-0.05289, -0.01162) at (-0.05289, -0.01162);
\coordinate (-0.13622, -0.03940) at (-0.13622, -0.03940);
\coordinate (-0.38622, -0.16440) at (-0.38622, -0.16440);
\coordinate (-1.00000, -1.99994) at (-1.00000, -1.99994);
\coordinate (-1.25000, -1.62494) at (-1.25000, -1.62494);
\coordinate (-1.33333, -1.40272) at (-1.33333, -1.40272);
\coordinate (-1.36458, -1.28553) at (-1.36458, -1.28553);
\coordinate (-1.37708, -1.22553) at (-1.37708, -1.22553);
\coordinate (-1.38229, -1.19515) at (-1.38229, -1.19515);
\coordinate (-1.38452, -1.17984) at (-1.38452, -1.17984);
\coordinate (-1.38550, -1.17215) at (-1.38550, -1.17215);
\coordinate (-1.38593, -1.16824) at (-1.38593, -1.16824);
\coordinate (-1.38613, -1.16629) at (-1.38613, -1.16629);
\coordinate (-1.38622, -1.16531) at (-1.38622, -1.16531);
\coordinate (-1.38622, -1.16440) at (-1.38622, -1.16440);
\fill[facet] (-1.38622, -1.16440) -- (-0.38622, -0.16440) -- (-0.13622, -0.03940) -- (-0.05289, -0.01162) -- (-0.02164, -0.00381) -- (-0.00914, -0.00131) -- (-0.00393, -0.00044) -- (-0.00169, -0.00012) -- (-0.00072, 0.00000) -- (0.00000, 0.00000) -- (0.00000, -1.99994) -- (-1.00000, -1.99994) -- (-1.25000, -1.62494) -- (-1.33333, -1.40272) -- (-1.36458, -1.28553) -- (-1.37708, -1.22553) -- (-1.38229, -1.19515) -- (-1.38452, -1.17984) -- (-1.38550, -1.17215) -- (-1.38593, -1.16824) -- (-1.38613, -1.16629) -- (-1.38622, -1.16531) -- cycle {};
\draw[edge] (0.00000, -1.99994) -- (0.00000, 0.00000);
\draw[edge] (0.00000, -1.99994) -- (-1.00000, -1.99994);
\draw[edge] (0.00000, 0.00000) -- (-0.00072, 0.00000);
\draw[edge] (-0.00072, 0.00000) -- (-0.00169, -0.00012);
\draw[edge] (-0.00169, -0.00012) -- (-0.00393, -0.00044);
\draw[edge] (-0.00393, -0.00044) -- (-0.00914, -0.00131);
\draw[edge] (-0.00914, -0.00131) -- (-0.02164, -0.00381);
\draw[edge] (-0.02164, -0.00381) -- (-0.05289, -0.01162);
\draw[edge] (-0.05289, -0.01162) -- (-0.13622, -0.03940);
\draw[edge] (-0.13622, -0.03940) -- (-0.38622, -0.16440);
\draw[edge] (-0.38622, -0.16440) -- (-1.38622, -1.16440);
\draw[edge] (-1.00000, -1.99994) -- (-1.25000, -1.62494);
\draw[edge] (-1.25000, -1.62494) -- (-1.33333, -1.40272);
\draw[edge] (-1.33333, -1.40272) -- (-1.36458, -1.28553);
\draw[edge] (-1.36458, -1.28553) -- (-1.37708, -1.22553);
\draw[edge] (-1.37708, -1.22553) -- (-1.38229, -1.19515);
\draw[edge] (-1.38229, -1.19515) -- (-1.38452, -1.17984);
\draw[edge] (-1.38452, -1.17984) -- (-1.38550, -1.17215);
\draw[edge] (-1.38550, -1.17215) -- (-1.38593, -1.16824);
\draw[edge] (-1.38593, -1.16824) -- (-1.38613, -1.16629);
\draw[edge] (-1.38613, -1.16629) -- (-1.38622, -1.16531);
\draw[edge] (-1.38622, -1.16531) -- (-1.38622, -1.16440);
\node[vertex] at (0.00000, 0.00000)     {};
\end{tikzpicture}
  \caption{An approximation to the triangle body $C$ obtained by weighted summation of the triangles $C^\tidx\ell$. The support set in direction $e_2$ consists of a single point marked black.}
  \end{subfigure}

\begin{center}
  \begin{subfigure}[b]{0.9\textwidth}
    \centering
    \begin{tikzpicture}%
	[scale=2.000000,
	back/.style={loosely dotted, thin},
	edge/.style={color=black!95!black, thick, dashed},
	facet/.style={fill=blue!95!black,fill opacity=0.800000},
	vertex/.style={inner sep=1pt,circle,draw=black!25!black,fill=black!75!black,thick,anchor=base}]
%
%
\coordinate (0.00000, -1.00000) at (0.00000, -1.00000);
\coordinate (0.00000, 0.00000) at (0.00000, 0.00000);
\coordinate (-1.00000, -1.00000) at (-1.00000, -1.00000);
\draw[edge] (0.00000, 0.00000) -- (-1.00000, -1.00000);
\node[vertex] at (0.00000, 0.00000)     {};
\node[vertex] at (-1.00000, -1.00000)     {};
\end{tikzpicture}\begin{tikzpicture}%
	[scale=2.000000,
	back/.style={loosely dotted, thin},
	edge/.style={color=black!95!black, thick, dashed},
	facet/.style={fill=blue!95!black,fill opacity=0.800000},
	vertex/.style={inner sep=1pt,circle,draw=black!25!black,fill=black!75!black,thick,anchor=base}]
%
%
\coordinate (0.00000, -2.00000) at (0.00000, -2.00000);
\coordinate (0.00000, 0.00000) at (0.00000, 0.00000);
\coordinate (-1.00000, -0.50000) at (-1.00000, -0.50000);
\draw[edge] (0.00000, 0.00000) -- (-1.00000, -0.50000);
\node[vertex] at (0.00000, 0.00000)     {};
\node[vertex] at (-1.00000, -0.50000)     {};
\end{tikzpicture}\begin{tikzpicture}%
	[scale=2.000000,
	back/.style={loosely dotted, thin},
	edge/.style={color=black!95!black, thick, dashed},
	facet/.style={fill=blue!95!black,fill opacity=0.800000},
	vertex/.style={inner sep=1pt,circle,draw=black!25!black,fill=black!75!black,thick,anchor=base}]
%
%
\coordinate (0.00000, -3.00000) at (0.00000, -3.00000);
\coordinate (0.00000, 0.00000) at (0.00000, 0.00000);
\coordinate (-1.00000, -0.33333) at (-1.00000, -0.33333);
\draw[edge] (0.00000, 0.00000) -- (-1.00000, -0.33333);
\node[vertex] at (0.00000, 0.00000)     {};
\node[vertex] at (-1.00000, -0.33333)     {};
\end{tikzpicture}\begin{tikzpicture}%
	[scale=2.000000,
	back/.style={loosely dotted, thin},
	edge/.style={color=black!95!black, thick, dashed},
	facet/.style={fill=blue!95!black,fill opacity=0.800000},
	vertex/.style={inner sep=1pt,circle,draw=black!25!black,fill=black!75!black,thick,anchor=base}]
%
%
\coordinate (0.00000, -4.00000) at (0.00000, -4.00000);
\coordinate (0.00000, 0.00000) at (0.00000, 0.00000);
\coordinate (-1.00000, -0.25000) at (-1.00000, -0.25000);
\draw[edge] (0.00000, 0.00000) -- (-1.00000, -0.25000);
\node[vertex] at (0.00000, 0.00000)     {};
\node[vertex] at (-1.00000, -0.25000)     {};
\end{tikzpicture}\begin{tikzpicture}%
	[scale=2.000000,
	back/.style={loosely dotted, thin},
	edge/.style={color=black!95!black, thick, dashed},
	facet/.style={fill=blue!95!black,fill opacity=0.800000},
	vertex/.style={inner sep=1pt,circle,draw=black!25!black,fill=black!75!black,thick,anchor=base}]
%
%
\coordinate (0.00000, -5.00000) at (0.00000, -5.00000);
\coordinate (0.00000, 0.00000) at (0.00000, 0.00000);
\coordinate (-1.00000, -0.20000) at (-1.00000, -0.20000);
\draw[edge] (0.00000, 0.00000) -- (-1.00000, -0.20000);
\node[vertex] at (0.00000, 0.00000)     {};
\node[vertex] at (-1.00000, -0.20000)     {};
\end{tikzpicture}\begin{tikzpicture}%
	[scale=2.000000,
	back/.style={loosely dotted, thin},
	edge/.style={color=black!95!black, thick, dashed},
	facet/.style={fill=blue!95!black,fill opacity=0.800000},
	vertex/.style={inner sep=1pt,circle,draw=black!25!black,fill=black!75!black,thick,anchor=base}]
%
%
\coordinate (0.00000, -6.00000) at (0.00000, -6.00000);
\coordinate (0.00000, 0.00000) at (0.00000, 0.00000);
\coordinate (-1.00000, -0.16667) at (-1.00000, -0.16667);
\draw[edge] (0.00000, 0.00000) -- (-1.00000, -0.16667);
\node[vertex] at (0.00000, 0.00000)     {};
\node[vertex] at (-1.00000, -0.16667)     {};
\end{tikzpicture}
    \caption{When the dashed faces of $C^\tidx\ell$ are scaled by a factor of $\ell$, they converge to a non-degenerate segment orthogonal to $e_2$.}
  \end{subfigure}
  \caption{The situation of Example \ref{ex:prune}}
  \label{fig:ex2}
  \end{center}
\end{figure}
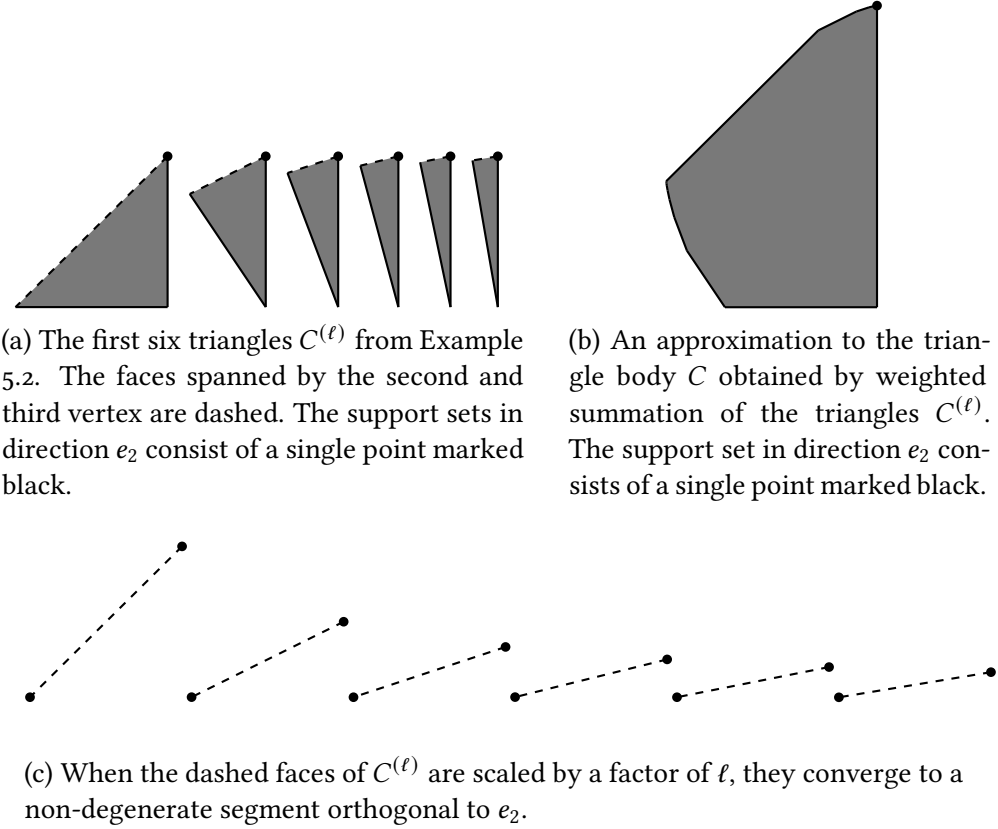

\begin{example}[Pruning]\label{ex:prune}{\rm
  Let $n=2$. Let again $(e_1,e_2)$ be the standard orthonormal basis of $\R^2$.
 Let
  \[
    C^\tidx\ell \coloneqq \conv\set*{-e_2,0, -\ell^{-1}e_1-\ell^{-2}e_2},\quad \ell\in\N
  ;\]
  see Figure \ref{fig:ex2} (a) and (b).
  Then we define $C$ to be the $3$-polyoid with generating measure
  \[
    \mu \coloneqq \sum_{\ell=1}^\infty 2^{-\ell} \delta_{C^\tidx\ell}
  .\]
  The sequence $(C^\tidx\ell)_\ell$ converges to the segment
  \[
    K \coloneqq \conv\set*{0,-e_2}
  ,\]
  so that $\supp\mu = \set*{K, C^\tidx1,C^\tidx 2, \ldots}$.
  This time, $\TS(C^\tidx\ell, e_2)$ for $\ell\in\N$ and $\TS(K, e_2)$ are all $0$-dimensional. However, there is no fixed $c > 0$ such that every $C^\tidx\ell$ has a $c$-cusp in direction $e_2$, so that $\TS(C, e_2)$ is nontrivial by Lemmas \ref{thm:tc4} and \ref{thm:tc5}.   Theorem \ref{thm:suppChar} predicts again that $e_2 \in \supp\Su(C)$.

  However, since $e_2 \notin \supp\Su(C^\tidx\ell)$ for all $\ell\in\N$ and $e_2 \notin \supp\Su(K)$, we cannot choose a \enquote{witness polytope} in $\supp\mu$ and repeat the argument from Example \ref{ex:simpleWitness}.

  The problem is this. In the previous example, the faces between the second and third vertex of $C^\tidx\ell$ converged to a one-dimensional face of the limit triangle $K$ with normal $e_2$. In the current example, however, these faces degenerate to a $0$-dimensional face. The only glimmer of hope is that the outer normals of these degenerating faces do still converge to $e_2$. If we could just scale up $C^\tidx\ell$ by a factor of $\ell$, the faces would not degenerate, but then we are confronted with the problem that $C^\tidx\ell$ is an unbounded sequence of convex bodies that does not converge to anything we might call a \enquote{witness polytope} anymore.

  On the other hand, by Lemma \ref{thm:facialStability} we find a neighborhood $U \subseteq \stdsph$ of $e_2$ such that for large enough $\ell$ and all $v \in U$,
  \[
    F(C^\tidx\ell, v) \subseteq \conv\set*{0, -\ell^{-1}e_1-\ell^{-2}e_2 } \eqqcolon F^\tidx\ell
  ;\]
  see Figure \ref{fig:ex2} (c) for an illustration.
  Therefore, for all Borel sets $V \subseteq U$,
  \[
    \tau(C^\tidx\ell, V) = \tau(F^\tidx\ell, V)
  .\]
  Now Lemma \ref{thm:maTau} implies that
  \[
    \Su(C^\tidx\ell)\mres{U} = \Su(F^\tidx\ell)\mres{U}
  .\]
  So if we can show that $e_2 \in \cl\bigcup_{\ell=1}^\infty\supp\Su(F^\tidx\ell)$, then $e_2 \in \cl\bigcup_{\ell=1}^\infty\supp\Su(C^\tidx\ell) = \supp\Su(C)$.
  Indeed, $\ell\cdot F^\tidx\ell\to F:=\conv\set*{0,-e_1}$, and therefore Lemma \ref{thm:suppEqOfTendsto} yields
  \[
    e_2 \in \supp\Su(F) \subseteq \cl\bigcup_{\ell=1}^\infty \supp\Su(\ell\cdot F^\tidx\ell) = \cl\bigcup_{\ell=1}^\infty \supp\Su(F^\tidx\ell)
  .\]
  In this example, we have leveraged that the $c$-cusps of $C^\tidx\ell$ in direction $e_2$ become more and more obtuse in the sense that $c > 0$ becomes smaller and smaller. This helped us find a sequence of faces, which unfortunately degenerated to a $0$-dimensional face in the limit. After \enquote{pruning} the sequence of triangles, i.e.\ removing some irrelevant vertices, we were able to scale up the polytopes in the sequence so that the sequence of faces converged to a $1$-dimensional face $F$, which we used as our \enquote{witness polytope} to prove that $e_2 \in \supp\Su(C)$.}
\end{example}

\begin{example}[Double pruning]{\rm
  We consider again $\R^2$ with the standard orthonormal basis $(e_1,e_2)$.
  Define for all $\ell\in\N$,
  \[
    v_1^\tidx\ell \coloneqq -e_2,\quad
    v_2^\tidx\ell \coloneqq 0,\quad
    v_3^\tidx\ell \coloneqq -\ell^{-1}e_1-\ell^{-1}e_2,\quad
    v_4^\tidx\ell \coloneqq -\ell^{-2}e_1-\ell^{-3}e_2,
  \]
  \[
    C^\tidx\ell \coloneqq \conv\set*{ v_1^\tidx\ell, v_2^\tidx\ell, v_3^\tidx\ell, v_4^\tidx\ell }
  .\]
  The vertices $v_2^\tidx\ell , v_3^\tidx\ell , v_4^\tidx\ell $ all converge to $0$, which is the unique element of the support set $F(\lim C^\tidx\ell, e_2)$. In analogy to the previous example, we remove $v_1^\tidx\ell $ and scale by a factor of $\ell$ to obtain a sequence of triangles
  \[
    D^\tidx\ell \coloneqq \conv\set*{ 0, -e_1-e_2,  -\ell^{-1}e_1-\ell^{-2}e_2}
  .\]
  Again, $F(\lim D^\tidx\ell, e_2)$ is a singleton and the vertices $\ell v_2^\tidx\ell, \ell v_4^\tidx\ell$ of $D^\tidx\ell$ converge to its unique element. Removing $\ell v_3^\tidx\ell$ and scaling by $\ell$ again, we get
  \[
    E^\tidx\ell \coloneqq \conv\set*{ 0, -e_1-\ell^{-1}e_2 }
  .\]
  Now, $F(\lim E^\tidx\ell, e_2) = \lim E^\tidx\ell$ is one-dimensional. Applying similar arguments as in the previous example, we conclude from $e_2 \in \cl\bigcup_{\ell=1}^\infty \supp\Su(E^\tidx\ell)$ that $e_2 \in \supp\Su(C)$.

  This example shows that the pruning procedure may have to be repeated several times.}
\end{example}

After these preparatory examples, we describe the general approach.

\begin{definition}
\label{def:prune}{\rm
  Let $\Qseq = (Q_\ell)_\ell$ be a bounded sequence of polytopes with a uniformly bounded number of vertices and $u \in \stdsph$. Let $k\in\N$ be the smallest number such that all polytopes in $\Qseq$ are $k$-topes.

  Choose an arbitrary sequence $\Vseq = (V_\ell)_\ell = ( (v_\ell^\tidx 1, \ldots, v_\ell^\tidx k) )_\ell$ of $k$-tuples of points in $\R^n$ such that
  \[
    Q_\ell = \conv\set*{ v_\ell^\tidx i \given i \in [k] } \quad \text{for all $\ell\in\N$}
  .\]
  Let $\Vseq' = (V_{\ell_s})_s$ be a convergent subsequence of $\Vseq$ and $Q\coloneqq \lim_{t\to\infty} Q_{\ell_t}$.
  Then we define a sequence $\prune(\Qseq, u) = (\prune(\Qseq, u, s))_s$ of polytopes
  \begin{align*}
    \prune(\Qseq, u, s) \coloneqq &\, c_s\left(\conv\set*{ v_{\ell_s}^\tidx i \given i \in [k], \lim_{t\to\infty} v_{\ell_t}^\tidx{i} \in F(Q, u) }- v_{\ell_s}^\tidx{i_0}\right)\\ \subseteq &\, c_s \left(Q_{\ell_s}-v_{\ell_s}^\tidx {i_0}\right)
  ,\end{align*}
where $i_0\in[k]$ is chosen such that $ \lim_{t\to\infty} v_{\ell_t}^\tidx{i_0} \in F(Q, u)$  and $c_s$ is the unique positive number such that $\diam\prune(\Qseq, u, s) =1$ if the convex hull by which $\prune(\Qseq, u, s)$ is defined is not a singleton, otherwise we set $c_s \coloneqq 1$, for $s\in\N$. Note that $0\in  \prune(\Qseq, u, s)$. We may  also pass to a subsequence of $\prune(\Qseq, u)$ and denote it in the same way; in any case, the sequence $\prune(\Qseq, u)$ is subject to various choices and not uniquely determined by $\Qseq$ and $u$. The polytopes in  $\prune(\Qseq, u)$ have diameter $1$ or are singletons and they contain $0$.

  If
  $$
  \lim_{t\to\infty} v_{\ell_t}^\tidx i \notin F(Q, u) \quad \text{for some }i\in[k],
  $$  then $k \ge 2$ and $\prune(\Qseq, u)$ consists of $(k-1)$-topes.
  After finitely many steps, the members of the sequence of sequences defined by
  \[
    \prune_0(\Qseq, u) \coloneqq \Qseq, \quad\prune_{m+1}(\Qseq, u) \coloneqq \prune(\prune_m(\Qseq, u), u) \quad \text{for all $m\in\N$}
  \]
  remain unchanged (if we do not pass to a subsequence) and become equal to some \enquote{fixpoint} sequence $\strip_*(\Qseq, u)$.}
\end{definition}

\begin{remark}\label{rem:subspace} {\rm
  If  $\strip_*(\Qseq, u)$ is obtained as described in Definition \ref{def:prune} and $Q^*:=\lim_{s\to\infty}\strip_*(\Qseq, u,s)$, then $0\in  Q^*\subset u^\perp$ and $\diam Q^*\in\{0,1\}$.}
\end{remark}

The next two lemmas prepare the proof of the crucial \nref{thm:pruningLemma}. The first is Lemma \ref{thm:realPruningLemma} which implies that at least locally pruning does not change the mixed area measures as far as their support is concerned. Lemma \ref{thm:sticky} then states a condition, which can be used to  ensure that the limit of a pruning sequence is non-degenerate.

\begin{lemma}[Pruning lemma]
  \makeatletter\def\@currentlabelname{Pruning Lemma}\makeatother\label{thm:realPruningLemma}
  Let $\Qseq = (Q_\ell)_\ell$ be a bounded sequence of polytopes in $\R^n$ with a uniform bound on the number of vertices, let $u \in\stdsph$ and $m \in \N_0$. Then there are an  $\stdsph$-open neighborhood $U \subseteq \stdsph$ of $u$, a subsequence $(Q_{\ell_s})_s$ and a sequence of positive numbers $(\lambda_s)_s$ such that for all but finitely many $s\in\N$ and for all $(n-2)$-tuples $\collc$ of  convex bodies in $\R^n$,
  \[
    \Su(Q_{\ell_s}, \collc)\mres U = \lambda_s \Su(\prune_m(\Qseq, u, s), \collc)\mres U
  .\]
  In particular, the statement is true if \ $\prune_m$ is replaced by $\prune_*$.
\end{lemma}
\begin{proof}
  The proof is by induction on $m \in \N_0$. If $m = 0$, the claim follows from $\Qseq = \prune_0(\Qseq, u)$. Now assume  $m \ge 1$ and that the claim is true for smaller $m$.

  Let $k\in\N$ be the smallest possible number such that $\Qseq$ consists of $k$-topes, just as in Definition \ref{def:prune}.
  Let $\Vseq = (V_\ell)_\ell = ( (v_\ell^\tidx 1, \ldots, v_\ell^\tidx k) )_\ell$ be a sequence of spanning points and $\Vseq' = (V_{\ell_s})_s$ a convergent subsequence (as in Definition \ref{def:prune}).
  Write
  \[
    v^\tidx i \coloneqq \lim_{s\to\infty} v_{\ell_s}^\tidx i \quad \text{for all $i \in [k]$}
  .\]
  We apply Lemma \ref{thm:facialStability} to $\lim_{s\to\infty} Q_{\ell_s} = \conv\set*{ v^\tidx i \given i \in [k] }$. Let
  \[
    I \coloneqq \set*{ i \in [k] \given v^\tidx i \in F\prn*{\lim_{s\to\infty} Q_{\ell_s}, u} }
  .\]
  Lemma \ref{thm:facialStability} shows that there is $\varepsilon\in (0,1) $ such that for all $w, x_1, \ldots, x_k$ with $d(u, w) < \varepsilon$ and $d(v^\tidx i, x_i) < \varepsilon$ ($i \in [k]$), the polytope $P \coloneqq \conv\set*{ x_i \given i \in [k] }$ satisfies
  \[
    F(P, w) \subseteq \conv\set*{ x_i \given i \in I }
  .\]
  In particular, there is an open neighborhood $U \subseteq \R^n\elminus0$ of $u$ such that for all but finitely many $s$ and for all $w \in U$,
  \[
    F(Q_{\ell_s}, w) -v_{\ell_s}^\tidx{i_0}\subseteq \conv\set*{ v_{\ell_s}^\tidx i \given i \in I } -v_{\ell_s}^\tidx{i_0}= c_s^{-1} \prune(\Qseq, u, s) \subseteq Q_{\ell_s}-v_{\ell_s}^\tidx{i_0}
  ,\]
  where $c_s$ is the positive factor in Definition \ref{def:prune}.
  It follows that
  \[F(\prune(\Qseq, u, s),w)=c_s\left(F(Q_{\ell_s},w)-v_{\ell_s}^\tidx{i_0}\right),
  \]
  and by Lemma \ref{thm:maTau} and the translation invariance of mixed area measures we get, for all but finitely many $s\in\N$ and for every $(n-2)$-tuple $\collc$ of  convex bodies,
  \[
    \Su(Q_{\ell_s}, \collc)\mres{(U \cap \stdsph)} = c_s^{-1} \Su(\prune(\Qseq, u, s), \collc)\mres{(U \cap \stdsph)}
  .\]

  Applying the inductive hypothesis for $m-1$ to $\prune(\Qseq, u)$, we obtain an $\stdsph$-open neighborhood $V \subseteq \stdsph$ of $u$, a subsequence $(\prune(\Qseq, u, s_t))_t$ and a sequence of positive numbers $(\mu_t)_t$ such that for all but finitely many $t \in \N$ and for every $(n-2)$-tuple $\collc$ of  convex bodies,
  \[
    \Su(\prune(\Qseq, u, s_t), \collc)\mres{V} = \mu_t \Su(\prune_m(\Qseq, u, t), \collc)\mres{V}
  .\]
  Now, $V \cap U \subseteq \stdsph$ is also an $\stdsph$-open neighborhood of $u$, and for all but finitely many $t \in \N$,
  \begin{align*}
    \Su(Q_{\ell_{s_t}}, \collc)\mres{(V \cap U)} &= c_{s_t}^{-1} \Su(\prune(\Qseq, u, s_t), \collc)\mres{(V\cap U)} \\
    &= c_{s_t}^{-1} \mu_t \Su(\prune_m(\Qseq, u, t), \collc)\mres{(V \cap U)}
  \end{align*}
  for every  $(n-2)$-tuple $\collc$ of convex bodies. This concludes the induction.

  Because there is $m\in\N_0$ such that $\prune_*(\Qseq, u) = \prune_m(\Qseq, u)$, the claim is also true for $\prune_*$.
\end{proof}

\begin{lemma}[Sticky vertices]\label{thm:sticky}
  Let $u \in\stdsph$. Let $\Qseq = (Q_\ell)_\ell$ be a bounded sequence of polytopes with a uniform bound on the number of vertices, satisfying the following property $ \mathfrak{P}(\Qseq,u)$: For all but finitely many $\ell\in\N$, there are distinct vertices $x_\ell, y_\ell$ of $Q_\ell$ such that $x_\ell \in F(Q_\ell, u)$ and $ \norm{x_\ell - y_\ell}^{-1}{\left<x_\ell - y_\ell, u\right>} \to 0$ as $\ell\to\infty$.

  Then  the sequences $\prune_m(\Qseq, u)$, $m\in\N$, can be chosen in such a way that the property $ \mathfrak{P}(\prune(\Qseq, u)_m,u) $ is satisfied for all $m\in\N$.
\end{lemma}
\begin{proof}
  Let $k\in\N$ be the smallest possible number such that $\Qseq$ consists of $k$-topes.
  It suffices to prove the claim for $\prune$, since the argument can then be iterated.

  Let $\Vseq = (V_\ell)_\ell = ( (v_\ell^\tidx 1, \ldots, v_\ell^\tidx k) )_\ell$ be a sequence of spanning points chosen as in Definition \ref{def:prune} which has $\Vseq' = (V_{\ell_s})_s$ as a   convergent subsequence. Moreover, the subsequence can be chosen such that there are distinct $i,j\in[k]$ with $x_{\ell_s}=v_{\ell_s}^\tidx i$ and $y_{\ell_s}=v_{\ell_s}^\tidx j$ for $s\in\N$. Then we set $Q:=\lim_{s\to\infty}Q_{\ell_s}$ and
  \[
    v^\tidx {i'} \coloneqq \lim_{s\to\infty} v_{\ell_s}^\tidx {i'} \quad \text{for $i' \in [k]$}
  .\]
It follows that
  \[ \left\langle v^\tidx i, u\right\rangle
   \leftarrow \left<v_{\ell_s}^\tidx i, u\right> = \left<x_{\ell_s}, u\right> = h_{Q_{\ell_s}}(u)\to h_Q(u),
  \]
as $s\to\infty$,
  hence  $\left<v^\tidx i, u\right> = h_Q(u)$ and    $v^\tidx i \in F(Q, u)$. Moreover,
  \[
   \left\langle v^\tidx j, u\right\rangle
   \leftarrow \left<v_{\ell_s}^\tidx j, u\right> = \left<y_{\ell_s}, u\right>
   = h_{Q_{\ell_s}}(u) + \left<y_{\ell_s} - x_{\ell_s}, u\right>\to h_Q(u),
  \]
 as $s\to\infty$, because of the assumption and since $\left(\left\|y_{\ell_s} - x_{\ell_s}\right\|\right)_s$ is bounded. Hence, we also have  $v^\tidx j \in F(Q, u)$.
  The construction of $\prune(\Qseq, u, s)$ then shows that $c_s (v_{\ell_s}^\tidx i-v_{\ell_s}^\tidx{i_0}) $ and $c_s (v_{\ell_s}^\tidx j-v_{\ell_s}^\tidx{i_0}) $ are distinct vertices of  $\prune(\Qseq, u, s)$ for all $s\in\N$, where $c_s$ is the positive scaling factor in Definition \ref{def:prune}.
  for  $s \in \N$.   In addition,
  $c_s (v_{\ell_s}^\tidx i-v_{\ell_s}^\tidx{i_0}) \in F(\prune(\Qseq, u, s),u)$ and
  \[
    \frac{\left<c_s (v_{\ell_s}^\tidx i-v_{\ell_s}^\tidx{i_0}) -c_s (v_{\ell_s}^\tidx j-v_{\ell_s}^\tidx{i_0}) , u\right>}{\norm{c_s (v_{\ell_s}^\tidx i-v_{\ell_s}^\tidx{i_0})  -c_s (v_{\ell_s}^\tidx j-v_{\ell_s}^\tidx{i_0}) }} = \frac{\left<x_{\ell_s} - y_{\ell_s}, u\right>}{\norm{x_{\ell_s} - y_{\ell_s}}} \to 0
  ,\]
  as $s\to\infty$. Thus $\prune(\Qseq, u)$ has the required property and the iteration can be continued.
\end{proof}

After these preparations, we state and prove the main auxiliary result in this section.

\begin{lemma}[Witness lemma]\label{thm:defaultPolytope}
  \makeatletter\def\@currentlabelname{Witness Lemma}\makeatother\label{thm:pruningLemma}
 Let $u\in\mathbb{S}^{n-1}$. Let $M \subset \R^n$ be a $k$-polyoid with generating measure $\mu$. If \ $\TS(M, u) \ne \{0\}$, then there is a $k$-tope $\Re(M,u) \subset u^\perp$ with $\{0\}\subset  \Re(M,u)$ (that is not a singleton) such that for every $(n-2)$-tuple $\collc$ of  convex bodies in $\R^n$,
  \[
    u \in \supp\Su(\Re(M,u), \collc) \quad \text{implies} \quad u \in \supp\Su(M, \collc)
  .\]
\end{lemma}
\begin{proof}
  For every $Q \in \supp\mu$, choose an arbitrary vertex $x_Q \in F(Q, u)$.

  If $\TS(M, u) \ne \{0\}$, then Corollary \ref{cor:cusp} shows that for every $c > 0$ there is some $P \in \supp\mu$ \emph{not} having a $c$-cusp in direction $u$. Hence we can find a sequence of $k$-topes $\Qseq:=(Q_\ell)_\ell$ in $\supp\mu$ and a sequence of  vertices $y_\ell\in Q_\ell$, $\ell\in\N$, such that $x_{Q_\ell}\neq y_\ell$ and
  \[
    \norm{y_\ell - x_{Q_\ell}}^{-1}  \left<y_\ell - x_{Q_\ell}, u\right>\to 0\quad \text{as }\ell\to\infty.
  \]
  By Lemma \ref{thm:sticky}, $\prune_*(\Qseq, u, \ell)$ has at least two distinct vertices, and by Definition \ref{def:prune} diameter $1$, for all but finitely many $\ell$.

  So $\prune_*( \Qseq, u)$, being a convergent sequence of $k$-topes, converges to a $k$-tope $\Re(M,u)\subset u^\perp$ of diameter $1$ with $0\in \Re(M,u)$ (see Remark \ref{rem:subspace}). In particular, $\Re(M,u)$ is not a singleton.

  By \nref{thm:realPruningLemma}, there is a sequence of positive numbers $(\lambda_s)_s$, a subsequence $(Q_{\ell_s})_s$ of $\Qseq$ and an $\mathbb{S}^{n-1}$-open neighborhood $U \subseteq \stdsph$ of $u$ such that for an arbitrary $(n-2)$-tuple $\collc$ of  convex bodies,
  \begin{align}\label{eq:jump}
    \Su(Q_{\ell_s}, \collc)\mres U = \lambda_s \Su(\prune_*( \Qseq, u, s ), \collc)\mres U
  .\end{align}

  Now assume that $u \in \supp\Su(\Re(M,u), \collc)$. Then by continuity of $\Su$ and Lemma \ref{thm:suppEqOfTendsto},
  \[
    u \in \cl\bigcup_{s=1}^\infty \supp\Su(\prune_*( \Qseq, u, s ), \collc)
  ,\]
  and because $U$ is a neighborhood of $u$, eq.~\eqref{eq:jump} and Theorem \ref{thm:suppInt} now imply
  \[
    u \in \cl\bigcup_{s=1}^\infty \supp\Su(Q_{\ell_s}, \collc) \subseteq \supp\Su(M, \collc)
  ,\]
which proves the assertion.
\end{proof}

\section{Switching}\label{sec:6}

In this section, we provide a lemma that will be needed in the proof of our main result to carry out the induction step. Recall the conventions and the notation concerning tuples of sets introduced in Section \ref{sec:2}. As usual, a linear subspace $R$ of some ambient vector space is said to be trivial if $R=\{0\}$.

\begin{lemma}[Switching lemma]\label{thm:critSwitching}
  Assume that $n \ge 2$ and $u \in \stdsph$. Let $\mathbfcal{T} = (T_1, \ldots, T_{n-1})$ and $\mathbfcal{R} = (R_1, \ldots, R_{n-1})$ be tuples of linear subspaces of \  $u^\perp$ such that $\mathbfcal{T}$ is semicritical and $R_i$ is nontrivial for all $i\in[n-1]$. Then there are index sets $\varnothing \ne I \subseteq J \subseteq [n-1]$ such that $\mathbfcal{R}_I$ spans an $\abs{I}$-dimensional subspace and $\mathbfcal{R}_J + \mathbfcal{T}_{J^\complement}$ is semicritical.
\end{lemma}
\begin{proof}
  Denote by $\mathbfcal{S} = (T_1 + R_1, \ldots, T_{n-1} + R_{n-1})$ the tuple of elementwise sums of $\mathbfcal{T}$ and $\mathbfcal{R}$.
  Choose $J \subseteq [n-1]$ inclusion-maximal such that
  \begin{align}\label{eq:jprop}
    \V\prn*{\mathbfcal{R}_J + \mathbfcal{S}_{J^\complement}} > 0
  .\end{align}
  Such $J$ exists since $\mathbfcal{S} = \mathbfcal{R}_\varnothing + \mathbfcal{S}_{[n-1]}$ is semicritical: $T_i \subseteq S_i$ for $i \in [n-1]$, $\mathbfcal{T}$ is semicritical by assumption and hence Lemma \ref{thm:critSimple} (6) implies the assertion.
  Because $J$ is maximal, even $\mathbfcal{R}_J + \mathbfcal{T}_{J^\complement}$ is semicritical: Repeatedly applying Lemma \ref{thm:critAdd}, we find a set $K \subseteq J^\complement$ such that $\mathbfcal{R}_{J\cup K} + \mathbfcal{T}_{J^\complement \setminus K}$ is semicritical. But then also $\mathbfcal{R}_{J \cup K} + \mathbfcal{S}_{J^\complement \setminus K}$ is semicritical, forcing $K = \varnothing$ because $J$ is inclusion-maximal.

  Furthermore, let $I \subseteq [n-1]$ be inclusion-minimal such that $\mathbfcal{R}_{I \cap J} + \mathbfcal{S}_{I \setminus J}$ is subcritical. Such $I$ exists because $\mathbfcal{R}_J + \mathbfcal{S}_{J^\complement}$ is subcritical as an $(n-1)$-tuple of $u^\perp$-subspaces (since  $n\ge 2$), showing that at least $[n-1]$ satisfies the desired property. Note that $I\neq\varnothing$, since by definition an empty tuple is not subcritical.

  Then $E \coloneqq \pspan \prn*{\mathbfcal{R}_{I \cap J} + \mathbfcal{S}_{I \setminus J}}$ is $\abs{I}$-dimensional: On the one hand, $\mathbfcal{R}_{I \cap J} + \mathbfcal{S}_{I \setminus J}$ is semicritical by Lemma \ref{thm:critSimple} (1). On the other hand, it is subcritical by the construction of $I$. If it spanned a higher-dimensional subspace, then this tuple would have to contain an even smaller subcritical set, contradicting the minimality of $I$.

  By Lemma \ref{thm:critRed} and relation \eqref{eq:jprop}, it follows that
  \begin{align}\label{eq:scproj}
   \V\prn*{\pi_{E^\perp}(\mathbfcal{R}_{J \setminus I} + \mathbfcal{S}_{J^\complement \setminus I})} > 0
  .\end{align}

  It remains to show that $I \subseteq J$.
  Assume for a contradiction that, without loss of generality, $1 \in I \setminus J$.

  Because $I$ was chosen inclusion-minimally such that $\mathbfcal{R}_{I \cap J} + \mathbfcal{S}_{I \setminus J}$ is subcritical, it follows that $\mathbfcal{R}_{I \cap J} + \mathbfcal{S}_{I \setminus (J \cup \set*{1})}$ is critical and as $R_1$ is nontrivial,  Lemma \ref{thm:critSimple} (7) implies that $\mathbfcal{R}_{I \cap (J \cup \set{1})} + \mathbfcal{S}_{I \setminus (J \cup \set{1})}$ is semicritical. Since $1\in I\setminus J$ and $R_1\subseteq S_1\subseteq E$, $\mathbfcal{R}_{I \cap (J \cup \set{1})} + \mathbfcal{S}_{I \setminus (J \cup \set{1})}$ spans a subspace of $E$ of dimension $\abs{I} = \dim E$, in other words, it also spans $E$.
  By Lemma \ref{thm:critRed} and relation \eqref{eq:scproj}, it follows that
  \[
    \V\prn*{\mathbfcal{R}_{J\cup\set{1}} + \mathbfcal{S}_{(J \cup\set{1})^\complement}} > 0,
  \]
 that is, $\mathbfcal{R}_{J\cup\set{1}} + \mathbfcal{S}_{(J \cup\set{1})^\complement}$ is semicritical, in contradiction to the maximality of $J$ as expressed by  relation \eqref{eq:jprop}. This proves that $I \subseteq J$.

  Finally, we get $\dim \mathbfcal{R}_{I}=|I|$,  since $I\subseteq J$ and $\dim E=|I|$.
\end{proof}

\section{Proof of the characterization theorem}\label{sec:7}

Now we are ready to confirm Theorem \ref{thm:suppChar} for smooth convex bodies and polyoids, for which eq.~\eqref{eq:MGM} should be recalled.

\begin{proof}

First, observe that it suffices to prove the claim for tuples that only contain polyoids (macroids). Consider the case that $\mathbfcal{M}$ does \emph{not} solely consist of polyoids (macroids) and let $\mathbfcal{M}'$ be the tuple obtained from $\mathbfcal{M}$ by replacing all smooth bodies by $\stdb$.
Clearly, $\mathbfcal{M}'$ consists of polyoids  (macroids), and the claim for $\mathbfcal{M}$ is equivalent to the claim for $\mathbfcal{M}'$ by the following argument.
All smooth convex bodies have the same, $(n-1)$-dimensional, touching spaces. Therefore, $\cl\ext \mathbfcal{M} = \cl\ext \mathbfcal{M}'$.
We can assume that the smooth and strictly convex body, contained in $\mathbfcal{M}$ by assumption, is the first one. By \cite[Lem.~7.6.15]{Schneider}, $\supp \Su(\mathbfcal{M}) = \supp\Su(\stdb, \mathbfcal{M}_{\backslash 1})$. Now \cite[Cor.~14.3]{SvH23+} shows that we can replace the remaining smooth bodies in $\mathbfcal{M}_{\backslash 1}$ by $\stdb$, and we obtain $\supp\Su(\mathbfcal{M}) = \supp\Su(\stdb, \mathbfcal{M}'_{\backslash 1}) = \supp\Su(\mathbfcal{M}')$. Hence, it suffices to prove the claim for tuples that only contain polyoids  (macroids), such as $\mathbfcal{M}'$.

It remains to show for tuples $\mathbfcal{M}$ of polyoids that
$$
  \supp \Su(\mathbfcal{M}) = \cl\ext \mathbfcal{M}.
$$
For this we prove two inclusions.

\medskip

\noindent
  \enquote{$\subseteq $}: For this part of the argument, we only need the weaker assumption that $\mathbfcal{M}$ is a tuple of macroids.
  By Theorem \ref{thm:suppInt} and Lemma \ref{thm:suppCharPoly}, we get
  \[
    \supp \Su(\mathbfcal{M}) = \cl\bigcup_\summandcolls \supp\Su(\collp) = \cl\bigcup_\summandcolls \ext \collp
  .\]
  So it remains to verify that
  \[
    \cl\bigcup_\summandcolls \ext(\collp)  \subseteq \cl\ext \mathbfcal{M}
  .\]
 Let $\collp = (P_1, \ldots, P_{n-1}) \in\prod_{i=1}^{n-1}\supp\mu_i$. We claim that for all $u \in \stdsph$ and $i \in [n-1]$,
  \begin{align}\label{eq:scp1}
    \TS(P_i, u) \subseteq \TS(M_i, u)
  ,\end{align}
  which would imply $\ext\collp \subseteq \ext \mathbfcal{M}$  by Lemma \ref{thm:critSimple} (6) and conclude the proof (for zonoids, compare \eqref{eq:scp1} with  \cite[Lem.~3.2]{Schneider1988}).

  Set $W \coloneqq \TS(M_i, u)^\perp$ and note that $u \in W$. Then by Lemma \ref{thm:tc3}, relation \eqref{eq:scp1} is equivalent to $\TS_W(\pi_W(P_i), u) \subseteq \TS_W(\pi_W(M_i), u) = \set*{0}$.
Now $\pi_W(P_i)$ is in the support of the generating measure of $\pi_W(M_i)$ by Lemma \ref{thm:projVertoid} (here we only need the inclusion which holds for general macroids). Together with $\TS_W(\pi_W(M_i), u) = \set*{0}$ and Lemmas \ref{thm:tc4} and \ref{thm:tc5}, this implies $\TS_W(\pi_W(P_i), u) = \set*{0}$ and relation \eqref{eq:scp1}.

\medskip

\noindent
  \enquote{$\supseteq $}:
  The proof of this inclusion is by induction on $n$. The case $n=1$ follows from Remark \ref{rem:empty} and the fact that the empty tuple is semicritical, rendering every $u \in S^0$ extreme.

 Assume  $n \ge 2$ and that the claim is true for smaller $n$. Let $u\in \ext \mathbfcal{M}$ be given.
  The linear subspaces $\TS(M_i, u)\subseteq u^\perp$, $i \in [n-1]$, form a semicritical tuple since $u\in \ext \mathbfcal{M}$, in particular $\TS(M_i, u)\neq\{0\}$. Then the linear subspaces $ \vspan{\Re(M_i,u)}\subseteq u^\perp$, which were defined in \nref{thm:pruningLemma}, are nontrivial and $\{0\}\subset  \Re(M_i,u)$ by Lemma \ref{thm:defaultPolytope}.

  Define
  \begin{align*}
    \mathbfcal{D} \coloneqq&\ (\vspan \Re(M_1,u), \ldots, \vspan \Re(M_{n-1},u))\\
     =&\ (\TS(\Re(M_1,u), u), \ldots, \TS(\Re(M_{n-1},u), u))
  ,\end{align*}
  where the equality follows from Lemma \ref{thm:touchingConePoly}, since $\Re(M_i,u)$ are polytopes with $0\in \Re(M_i,u)\subseteq u^\perp$ so that $F(\Re(M_i,u),u)=\Re(M_i,u)$, for $i\in [n-1]$.
  According to Lemma \ref{thm:critSwitching}, there are index sets $\varnothing \ne I \subseteq J \subseteq [n-1]$ such that $\mathbfcal{D}_I$ spans an $\abs{I}$-dimensional subspace $E$ and $\mathbfcal{D}_J + \TS(\mathbfcal{M}_{J^\complement}, u)$ is semicritical.

  We now interpret the $k$-tope $\Re(M_i,u)$, where $i \in J$, as a $k$-polyoid with generating Dirac measure $\delta_{\Re(M_i,u)}$ and define
  \[
    \mathbfcal{M}' \coloneqq (M_1', \ldots, M_{n-1}'), \quad \text{where } M_i' \coloneqq
    \begin{cases}
      \Re(M_i,u),&   i \in J, \\
      M_i,&  i \notin J.
    \end{cases}
  \]
  It now suffices to prove that $u \in \supp\Su(\mathbfcal{M}')$: Using that $\Re(M_j,u) = M'_j$ ($j \in J$), repeated applications of \nref{thm:pruningLemma} show that if $u \in \supp\Su(\mathbfcal{M}')$, then $u \in \supp\Su(\mathbfcal{M})$.

  Clearly, $\mathbfcal{M}'$ is also a tuple of $k$-polyoids and $\TS( \mathbfcal{M}', u)$ is semicritical because $\mathbfcal{D}_J + \TS(\mathbfcal{M}_{J^\complement}, u)$ is a semicritical permutation.
  Furthermore, all spaces in $\TS(\mathbfcal{M}'_I, u) = \mathbfcal{D}_I$ are subspaces of $E$. Lemmas \ref{thm:critRed} and \ref{thm:tc3} now imply that $\pi_{E^\perp}\TS(\mathbfcal{M}'_{I^\complement}, u) = \TS_{E^\perp}(\pi_{E^\perp}\mathbfcal{M}'_{I^\complement}, u)$ is also a semicritical tuple, that is, we have $u\in \ext \pi_{E^\perp}(\mathbfcal{M}')_{I^\complement}$.

  There is an inner product space isomorphism $E^\perp \cong \R^{\dim E^\perp}$. Using this isomorphism and $\dim E^\perp = n - \abs{I} \in [1, n-1]$, we can apply the inductive hypothesis to $u \in E^\perp$ and the tuple $\pi_{E^\perp}(\mathbfcal{M}')_{I^\complement}$ of $k$-polyoids in $E^\perp$ and thus conclude from $u\in \ext \pi_{E^\perp}(\mathbfcal{M}')_{I^\complement}$ that
  \begin{align}\label{eq:mp1}
    u \in \supp\Su_{E^\perp}(\pi_{E^\perp}(\mathbfcal{M}')_{I^\complement})
  .\end{align}

  On the other hand, $\mathbfcal{M}_I'$ consists of $k$-topes in $E$. So Lemma \ref{thm:maRed} yields
  \begin{align}\label{eq:mp2}
    \begin{pmatrix}
      n-1 \\ \abs{I}
    \end{pmatrix}
    \Su(\mathbfcal{M}') = \V(\mathbfcal{M}'_I) \cdot \Su_{E^\perp}'(\pi_{E^\perp}(\mathbfcal{M}'_{I^\complement}))
  .\end{align}
  Because $(\vspan M'_i)_{i\in I} = \mathbfcal{D}_I$ is a subtuple of the semicritical tuple $\mathbfcal{D}_J + \TS(\mathbfcal{M}_{J^\complement}, u)$, since $I\subseteq J$, it follows from Lemma \ref{thm:mvVanish} that $\V(\mathbfcal{M}'_I) > 0$ and we conclude with relations \eqref{eq:mp1} and \eqref{eq:mp2} that
  \[
    u \in \supp\Su_{E^\perp}'(\pi_{E^\perp}(\mathbfcal{M}')_{I^\complement}) = \supp\Su(\mathbfcal{M}')
  \]
  and, as noted previously, therefore $u \in \supp\Su(\mathbfcal{M})$.
\end{proof}

Finally, we note the following more general result which is implied by the preceding proof.

\begin{proposition}
  \label{prop:suppChar}
  Let $\collc = (C_1, \ldots, C_{n-1})$ be an $(n-1)$-tuple of macroids (or smooth convex bodies provided at least one of the bodies $C_i$ is smooth and strictly convex) in $\R^n$. Then
  \begin{equation}\label{eq:key2b}
 \supp\Su(\collc,\cdot) \subseteq \cl\cextdirs
  .\end{equation}
\end{proposition}

\bigskip

\noindent
\textbf{Acknowledgements.}
D. Hug was supported by DFG research grant HU 1874/5-1 (SPP 2265).

\bigskip

\vspace{2cm}

\noindent
Authors' addresses:

\bigskip

\noindent
Daniel Hug, Karlsruhe Institute of Technology (KIT), Institute of Stochastics, D-76128 Karlsruhe, Germany, daniel.hug@kit.edu

\medskip

\noindent
Paul Reichert, Karlsruhe, Germany, math@paulr.de

\end{document}